\documentclass[12pt,a4paper]{amsart}
\usepackage{amssymb,amsmath}

\usepackage{times}

\numberwithin{equation}{section}
\newtheorem{theorem}{Theorem}[section]
\newtheorem{lemma}[theorem]{Lemma}
\newtheorem{proposition}[theorem]{Proposition}
\newtheorem{corollary}[theorem]{Corollary}

\theoremstyle{definition}
\newtheorem{definition}[theorem]{Definition}
\newtheorem{example}[theorem]{Example}
\newtheorem{remark}[theorem]{Remark}

\newcommand\Supp{\operatorname{Supp}}
\newcommand\Ass{\operatorname{Ass}}
\newcommand\Assh{\operatorname{Assh}}
\newcommand\Ann{\operatorname{Ann}}
\newcommand\Min{\operatorname{Min}}
\newcommand\Tor{\operatorname{Tor}}
\newcommand\Hom{\operatorname{Hom}}

\newcommand\Ext{\operatorname{Ext}}
\newcommand\Rad{\operatorname{Rad}}
\newcommand\Ker{\operatorname{Ker}}
\newcommand\Coker{\operatorname{Coker}}
\newcommand\im{\operatorname{Im}}

\newcommand{\xx}{\underline x}

\newcommand\grade{\operatorname{grade}}
\newcommand\fgrade{\operatorname{fgrade}}
\newcommand\cd{\operatorname{cd}}
\newcommand\depth{\operatorname{depth}}
\newcommand\height{\operatorname{height}}
\newcommand\Spec{\operatorname{Spec}}
\newcommand{\gam}{\Gamma_{\mathfrak a}}
\newcommand{\Rgam}{{\rm R} \Gamma_{\mathfrak a}}
\newcommand{\lam}{\Lambda^{\mathfrak a}}

\newcommand{\Cech}{{\Check {C}}_{\xx}}
\newcommand{\HH}{H_{\mathfrak m}}
\newcommand{\qism}{\stackrel{\sim}{\longrightarrow}}

\begin{document}

\author[P.~Schenzel]{Peter Schenzel}
\title[Formal cohomology]
{On formal local cohomology and connectedness}
\address{Martin-Luther-Universit\"at Halle-Wittenberg,
Institut f\"ur Informatik,
D --- 06 099 Halle (Saale), Germany}
\email{schenzel@mathematik.uni-halle.de}

\begin{abstract} Let $\mathfrak a$ denote an ideal of a local
ring $(R, \mathfrak m).$  Let $M$ be a finitely generated
$R$-module. There is a systematic study of the formal cohomology
modules $\varprojlim \HH^i(M/\mathfrak a^nM), i \in \mathbb Z.$ We
analyze their $R$-module structure, the upper and lower vanishing
and non-vanishing in terms of intrinsic data of $M,$ and its
functorial behavior. These cohomology modules occur in relation to
the formal completion of the punctured spectrum $\Spec R \setminus
V(\mathfrak m).$

As a new cohomological data there is a description on the formal
grade $\fgrade( \mathfrak a, M)$ defined as the minimal
non-vanishing of the formal cohomology modules. There are various
exact sequences concerning the formal cohomology modules. Among
them a Mayer-Vietoris sequence for two ideals. It applies to new
connectedness results. There are also relations to local
cohomological dimensions.
\end{abstract}

\maketitle

\section{Introduction} Let $\mathfrak a$ denote an ideal of a local ring
$(R, \mathfrak m).$ For a finitely generated $R$-module $M$ let $H^i_{\mathfrak a}(M),
 i \in \mathbb N,$ denote the local cohomology module of $M$ with respect to
$\mathfrak a$ (cf. \cite{aG} for the basic definitions). There are
the following integers related to these local cohomology modules
\begin{equation*}
\grade ( \mathfrak a, M)  =  \inf \{i \in \mathbb Z :
H^i_{\mathfrak a}(M) \not= 0\}, \; \cd (\mathfrak a, M)  = \sup
\{i \in \mathbb Z : H^i_{\mathfrak a}(M) \not= 0\},
\end{equation*}
called the grade (resp. the cohomological dimension) of $M$ with
respect $\mathfrak a$ (cf. Section 2.2). In general we have the
bounds $\height _M \mathfrak a \leq \cd(\mathfrak a, M) \leq \dim
M.$ In the case of $\mathfrak m$ the maximal ideal it follows that
$\grade (\mathfrak m, M) = \depth M$ and $\cd(\mathfrak m, M) =
\dim M.$

Here we consider the asymptotic behavior of the family of local
cohomology modules $\{\HH^i(M/\mathfrak a^n M)\}_{n \in \mathbb
N}$ for an integer $i \in \mathbb Z.$ By the natural homomorphisms
these families form a projective system. Their projective limit
$\varprojlim \HH^i(M/\mathfrak a^n M)$ is called the $i$-th formal
local cohomology of $M$ with respect to $\mathfrak a.$ Not so much
is known about these modules. In the case of a regular local ring
they have been studied by Peskine and Szpiro (cf. \cite[Chapter
III]{PS}) in relation to the vanishing of local cohomology
modules. Another kind of investigations about formal cohomology
has been done by Faltings (cf. \cite{gF}).

Moreover, $\varprojlim \HH^i(M/\mathfrak a^n M)$ occurs as the
$i$-th cohomology module of the $\mathfrak a$-adic completion of
the \v{C}ech complex $\Cech \otimes M$ (cf. Section 3), where
$\xx$ denotes a system of elements of $R$ such that $\Rad \xx R =
\mathfrak m.$

The main subject of the paper is a systematic study of the formal
local cohomology modules. Above all we are interested in the first
resp. last non-vanishing of the formal cohomology. As an easy result of
this type the following result is proved:

\begin{theorem} \label{1.1} Let $\mathfrak a$ denote an ideal of a local
ring $(R, \mathfrak m).$ Then
\[
\dim M/\mathfrak a M = \sup \{i \in \mathbb Z : \varprojlim
\HH^i(M/\mathfrak a^nM) \not= 0\}
\]
for a finitely generated $R$-module $M.$
\end{theorem}

The description of $\inf \{i \in \mathbb Z : \varprojlim
\HH^i(M/\mathfrak a^nM) \not= 0\}$ in terms of intrinsic data
seems to be not obvious. Following the intention of Peskine and
Szpiro (cf. \cite[Chapter III]{PS}) we define the formal grade as
\[
\fgrade (\mathfrak a, M) = \inf \{i \in \mathbb Z : \varprojlim
\HH^i(M/\mathfrak a^nM) \not= 0\}
\]
for an ideal $\mathfrak a$ and a finitely generated $R$-module $M.$
Since the formal cohomology does not change by passing to the
completion of $R$ (cf. \ref{3.3}) we may assume -- without loss of generality --
the existence of a
dualizing complex $D^{\cdot}_R$ for $R$ . So
we may express the formal cohomology in terms of the local cohomology
of the dualizing complex.

\begin{theorem} \label{1.2} Let $(R, \mathfrak m)$ denote a local ring possessing a
(normalized) dualizing complex $D^{\cdot}_R$. Let $\mathfrak a$ denote an ideal of $R.$ For a
finitely generated $R$-module $M$ it follows
\begin{itemize}
\item[(a)] $\varprojlim \HH^i(M/\mathfrak a^nM) \simeq
\Hom_R(H^{-i}_{\mathfrak a}(\Hom_R(M, D^{\cdot}_R)),E),$ for all
$i \in \mathbb Z,$
\item[(b)] $\fgrade (\mathfrak a, M) = \inf \{
i - \cd (\mathfrak a, K^i(M)) : i = 0, \ldots , \dim M\}.$
\end{itemize}
Here $K^i(M) = H^{-i}(\Hom(M, D^{\cdot}_R)), i = 0,\ldots, \dim M,$ denotes
the $i$-th module of deficiency (cf. Section 2.3).
\end{theorem}

Another result concerns the vanishing of the
formal cohomology $\varprojlim \HH^i(M/\mathfrak a^nM)$ and the dimension
of the associated prime ideals of the underlying module.

\begin{theorem} \label{1.3} Let $\mathfrak a$ denote an ideal of a
local ring $(R, \mathfrak m).$ Let $M$ be a finitely generated $R$-module. Then
\begin{itemize}
\item[(a)] $\fgrade(\mathfrak a, M) \leq \dim M - \cd(\mathfrak a,
M),$

\item[(b)] $\fgrade (\mathfrak a, M) \leq \dim \hat{R}/ \mathfrak
p - \cd(\mathfrak a \hat{R}, \hat{R}/\mathfrak p) \quad \mbox{ for
all } \quad \mathfrak p \in \Ass \hat{M},$
\end{itemize}
where $\hat {M}$ denotes the $\mathfrak m$-adic completion of $M.$
\end{theorem}

In Section 5 there is a Mayer-Vietoris sequence for the formal
cohomology, analogous to the corresponding sequence for the local
cohomology. As in the case of local cohomology this applies to
connectedness results of certain subsets of $\Spec R.$ To this end
let $c(R/\mathfrak c)$ denote the connectedness dimension of
$V(\mathfrak c)$ for an ideal $\mathfrak c$ (cf. \ref{5.9}).

\begin{theorem} \label{1.4} Let $\mathfrak a$ be an ideal of $(R,\mathfrak m).$
For a finitely generated $R$-module M there are the estimates:
\begin{itemize}
\item[(a)] $\fgrade (\mathfrak a, M)  - 1 \leq
 c(\hat{R}/(\mathfrak a \hat{R}, \mathfrak p))
\quad \mbox{ for all } \mathfrak p \in \Ass \hat{M}.$

\item[(b)] Assume that $\Ass \hat{M} = \Assh \hat{M}$ and
$\HH^d(R/\Ann M), d = \dim M,$ is an indecomposable $R$-module.
 Then $\fgrade (\mathfrak a, M) -1 \leq c(\hat{M}/\mathfrak a \hat{M}).$
\end{itemize}
\end{theorem}

In particular, when $\varprojlim \HH^i(M/\mathfrak a^ nM) = 0$ for $i = 0,1,$ then
$V(\mathfrak a \hat{R}, \mathfrak p) \setminus V(\hat{\mathfrak m})$ is connected
for all $\mathfrak p \in \Ass \hat{M}.$

In Section 2 of the paper we start with some preliminaries about notation,
local cohomology, dualizing complexes, and commutative algebra. Section 3 is devoted to the
definitions and basic results about formal cohomology, its relation to duality,
as well as exact sequences for various situations. In Section 4 there are vanishing
and non-vanishing results about formal cohomology. This Section contains also the
results about the formal grade. In Section 5 there is the Mayer-Vietoris sequence
for formal cohomology and the connectedness properties. In addition there are
also results about the connectedness and the local cohomology.

\section{Preliminary Results}
\subsection{Notation}
In the present paper $(R, \mathfrak m, k)$ denotes a local
Noetherian ring with its residue field $k = R/\mathfrak m.$ In the
following let $\mathfrak a, \mathfrak b, \ldots$ denote ideals of
$R.$ Let $M$ be an $R$-module. By $X: \ldots \rightarrow X^n
\stackrel{d_X^n}{\rightarrow} X^{n+1} \rightarrow \dots$ we denote
a complex of $R$-modules.

Let $\xx = x_1,\ldots,x_n$ be a sequence of elements of $R.$ Then
$K_{\cdot}(\xx; X)$ and $K^{\cdot}(\xx;X)$ are the Koszul
complexes of $X$ with respect to $\xx$ (cf. \cite{pS3} for the
definition of Koszul complexes and basic facts about homological
algebra).

For an arbitrary $R$-complex $X$ there is a complex $I$ of
injective $R$-modules (resp. a complex $F$ of flat $R$-modules)
and a quasi-isomorphism $X \qism I$ (resp. $F \qism X$) (cf.
\cite{nS} or \cite{AF} for the construction). We call $I$ (resp.
$F$) an injective (resp. a flat) resolution of $X.$

For an $R$-complex $X$ and an integer $m \in \mathbb Z$ define the
shifted complex $X[m]$ by $X[m]^n = X^{m+n}, n \in \mathbb Z,$ and
$d_{X[m]} = (-1)^m d_X,$ where $d$ denotes the boundary map.

\subsection{Local cohomology} Let $\xx = x_1,\ldots,x_n$ be a
system of elements of the ring $R$ and let $\mathfrak a =
(x_1,\ldots,x_n)R$ the ideal generated by these elements. The
local cohomology $\Rgam (X)$ of $X$ with respect to $\mathfrak a$
in the derived category  is defined by $\gam (I),$ where $X \qism
I$ denotes the injective resolution (cf. \cite{rH} resp.
\cite{hF}). For an integer $i \in \mathbb Z$ define
$H^i_{\mathfrak a}(X) = H^i(\gam (I)).$ Note that up to
isomorphisms it is independent on $I.$

Moreover let $\Cech$ denote the \v{C}ech complex with respect to
$\xx$ (cf. \cite{pS1} or \cite{pS2}). Then there is a canonical
isomorphism $\gam (I) \simeq \Cech \otimes I$ for a complex of
injective $R$-modules $I$ (cf. \cite{pS2}). Because $\Cech$ is a
bounded $R$-complex of flat $R$-modules it induces the following
isomorphism $\Cech \otimes X \simeq \Cech \otimes I.$ That is, the
local cohomology $H^i_{\mathfrak a}(X), i \in \mathbb Z,$ may be
computed as the cohomology $H^i(\Cech \otimes X).$

For a finitely generated $R$-module $M$ there is the following
characterization
\[
\grade ( \mathfrak a, M) = \inf \{i \in \mathbb Z : H^i_{\mathfrak
a}(M) \not= 0\}
\]
for the $\grade ( \mathfrak a,M)$ of the $R$-module $M$ with
respect to the ideal $\mathfrak a.$ For the supremum of the
non-vanishing there is the following definition
\[
\cd (\mathfrak a, M) = \sup \{i \in \mathbb Z : H^i_{\mathfrak a}(M)
\not= 0\},
\]
where $\cd ( \mathfrak a,M)$ is called the cohomological dimension
of $M$ with respect to $\mathfrak a.$ Recall that $\cd ( \mathfrak
a ,M)\leq \dim_R M$ with equality in the case $\Rad \mathfrak a =
\mathfrak m$ (cf. \cite{aG}). Moreover $\height_M \mathfrak a \leq
\cd( \mathfrak a, M),$ where $\height_M \mathfrak a = \height
(\mathfrak a, \Ann_R M)/\Ann_R M.$ In general it is a difficult
problem to calculate the cohomological dimension $\cd(\mathfrak
a,R)$ of an ideal.

We need here another preliminary result about cohomological
dimensions. It was invented by Divaani-Aazar, Naghipour and Tousi
(cf. \cite{DNT}). For sake of completeness we include a proof.

\begin{lemma} \label{5.3} Let $\mathfrak a$ denote an ideal of a local
ring $(R, \mathfrak m).$ Let $M, N$ be two finitely generated $R$-modules
such that $\Supp N \subseteq \Supp M.$ Then $\cd (\mathfrak a, N) \leq
\cd (\mathfrak a,  M).$
\end{lemma}

\begin{proof} It will be enough to show that $H^i_{\mathfrak a}(N) = 0$
for all integers $\cd (\mathfrak a, M) < i \leq \dim M + 1.$ The proof
will be shown by an descending induction on $i.$

First note that the claim is true for $i = \dim M +1.$ (cf.
\cite{aG}). Now let $i \leq \dim M.$ We proceed by a trick
invented by Delfino and Marley (cf. the proof of \cite[Proposition
1]{DM}). By the assumption we have $\Supp N \subseteq \Supp M,$
and therefore $\Rad \Ann_R N \supseteq \mathfrak c,$ where
$\mathfrak c = \Ann _RM.$ Whence there is an $n \in \mathbb N$
such that $\mathfrak c^ n N = 0.$ Thus $N$ possesses a filtration
\[
0 = \mathfrak c^ n N \subset \mathfrak c^{n-1} N \subset \ldots \subset
\mathfrak c N \subset N,
\]
such that $\mathfrak c^{i-1}N/\mathfrak c^iN, i = 1, \ldots n,$ is a
finitely generated $R/\mathfrak c$-module.

By Gruson's theorem (cf. \cite[Theorem 4.1]{wV}) a finitely
generated $R/\mathfrak c$-module $T$ admits a filtration
\[
0 = T_0 \subset T_1 \subset \ldots \subset T_k = T,
\]
such that $T_j/T_{j-1}, j = 1,\ldots,k,$ is a homomorphic image of finitely
many copies of $M.$

We prove now the vanishing of $H^ i_{\mathfrak a}(T).$ By using
short exact sequences and induction on $k$ it suffices to prove
the case when $k = 1.$ Thus, there is an exact sequence
\[
0 \to K \to M^m \to T \to 0
\]
for some positive integer $m.$ It induces an exact sequence
\[
\ldots \to H^i_{\mathfrak a}(K) \to  H^i_{\mathfrak a}(M)^m \to H^i_{\mathfrak a}(T)
\to H^{i+1}_{\mathfrak a}(K) \to \ldots .
\]
By the inductive hypothesis $H^{i+1}_{\mathfrak a}(K) = 0,$ so that $H^i_{\mathfrak a}(T) = 0.$

Finally we prove that $H^i_{\mathfrak a}(N) = 0.$ By the use of short exact sequences
and induction on $n,$ it suffices to prove the case when $n = 1,$ which is obviously true
by the aid of the previous argument.
\end{proof}

As a corollary of the previous Lemma \ref{5.3} it follows that the cohomological dimension
of a finitely generated $R$-module $M$ is determined by the cohomological
dimension of its minimal associated prime ideals. To this end let $\Min M$ denote the
minimal elements of $\Supp M,$ where $M$ denotes an $R$-module.

\begin{corollary} \label{5.4} Let $M$ be a finitely generated $R$-module. Then
\[
\cd (\mathfrak a, M) = \cd (\mathfrak a, R/\Ann_R M) = \max \{ \cd
(\mathfrak a, R/\mathfrak p) : \mathfrak p \in \Min M\}
\]
for any ideal $\mathfrak a$ of $R.$
\end{corollary}

\begin{proof} The fist equality is clear because of $V(\Ann_R M) = \Supp_RM$
(cf. \ref{5.3}). For the proof of the second define $N =
\oplus_{\mathfrak p \in \Min M} R/\mathfrak p.$ Then it follows
that
\[
\cd (\mathfrak a, N )= \max \{ \cd (\mathfrak a, R/\mathfrak p) :
\mathfrak p \in \Min M\}.
\]
Remember that the local cohomology commutes with direct sums.
Furthermore we have $\Supp M = \Supp N.$ So the statement is a
consequence of Lemma \ref{5.3}.
\end{proof}

As another preliminary result we need the behavior of the
cohomological dimension of an $R$-module with respect to an ideal
$\mathfrak a$ by passing to $(\mathfrak a, xR).$

\begin{lemma} \label{2.x} Let $\mathfrak a$ denote an ideal of a local ring
$(R, \mathfrak m).$ Let $M$ be a finitely generated $R$-module.
Then
\[
\cd( (\mathfrak a, xR), M) \leq \cd(\mathfrak a, M) + 1
\]
for any element $x \in \mathfrak m.$
\end{lemma}

\begin{proof} With the notation of the lemma there is the short
exact sequence
\[
0 \to H^1_{xR}(H^i_{\mathfrak a}(M)) \to H^{i+1}_{(\mathfrak a ,
xR)}(M) \to H^0_{xR}(H^{i+1}_{\mathfrak a}(M)) \to 0
\]
for all $i \in \mathbb Z$ (cf. for instance \cite[Corollary
3.5]{pS2}). Now put $c = \cd (\mathfrak a, M).$ Then by the
definition of the cohomological dimension the short exact sequence
implies that $H^{i+1}_{(\mathfrak a, xR)}(M) = 0$ for all $i > c.$
In other words $\cd ((\mathfrak a, xR), M) \leq c +1,$ which
finishes the proof.
\end{proof}

\subsection{Dualizing complexes}
In this subsection let $(R, \mathfrak m)$ denote a local ring
possessing a dualizing complex $D^{\cdot}_R.$ That is a bounded
complex of injective $R$-modules whose cohomology modules
$H^i(D^{\cdot}_R), i \in \mathbb Z,$ are finitely generated
$R$-modules. We refer to \cite[Chapter V, \S 2]{rH} or to
\cite[1.2]{pS1} for basic results about dualizing complexes.

By the result of T.~Kawasaki (cf. \cite{tK}) $R$ possesses a
dualizing complex if and only if $R$ is the factor ring of a
Gorenstein ring.

Note that the natural homomorphism of complexes
\[
M \to \Hom_R(\Hom_R(M, D^{\cdot}_R), D^{\cdot}_R)
\]
induces an isomorphism in cohomology for any finitely generated
$R$-module $M.$ Moreover there is an integer $l \in \mathbb Z$
such that
\[
\Hom_R(k, D^{\cdot}_R) \simeq k[l],
\]
where $k = R/\mathfrak m$ denotes the residue field of $R.$ As
follows by a shifting we may always assume without loss of
generality assume that $l = 0.$ Then the dualizing complex
$D^{\cdot}_R$ is called normalized. In the following let us always
assume that a dualizing complex is normalized.

Then a dualizing complex has the following structure
\[
D^{-i}_R \simeq \oplus_{\mathfrak p \in \Spec R, \dim R/\mathfrak
p =i} E_R(R/\mathfrak p),
\]
where $E_R(R/\mathfrak p)$ denotes the injective hull of
$R/\mathfrak p$ as $R$-module. Therefore $D^i_R = 0$ for $i < -
\dim R$ and $i > 0.$

\begin{proposition} \label{2.1} Let $(R, \mathfrak m)$ denote a
local ring with the dualizing complex $D^{\cdot}_R.$
\begin{itemize}
\item[(a)] $D^{\cdot}_R \otimes R_{\mathfrak p} \simeq
D^{\cdot}_{R_{\mathfrak p}}[\dim R/\mathfrak p]$ \; for \; $\mathfrak p
\in \Spec R.$
\item[(b)] {\rm (Local duality)} There is a
canonical isomorphism
\[
H^i_{\mathfrak m}(M) \simeq \Hom_R (H^{-i}(\Hom_R(M,
D^{\cdot}_R)), E), \quad E = E_R(R/\mathfrak m),
\]
for a finitely generated $R$-module $M$ and all $i \in \mathbb Z.$
\end{itemize}
\end{proposition}

The proof is well-known (cf. \cite{rH} resp. \cite{pS1}). For a
certain application remember the definition of the modules of
deficiencies of an $R$-module $M$ (cf. \cite[Section 1.2]{pS1}).

\begin{definition} \label{2.z} Let $M$ denote a finitely generated
$R$-module and $d = \dim M.$ For an integer $i \in \mathbb Z$ define
$$
K^i(M) := H^{-i}(\Hom_R(M, D^{\cdot}_R)).
$$
The module $K(M) := K^d(M)$ is called the canonical module of $M.$
For $i \not= d$ the modules $K^i(M)$ are called the modules of
deficiency of $M.$ Note that $K^i(M) = 0$ for all $i < 0$ or $i > d.$
\end{definition}

By the local duality theorem there are the canonical
isomorphisms
$$
\HH^i(M) \simeq \Hom_R(K^i(M), E), i \in \mathbb Z,
$$
where $E = E_R(R/\mathfrak m)$ denotes the injective hull of the
residue field. Remember that all of the $K^i(M), i \in \mathbb Z,$
are finitely generated $R$-modules. Moreover $M$ is a
Cohen-Macaulay module if and only if $K^i(M) = 0$ for all $i \not=
d.$ Whence the modules of deficiencies of $M$ measure the
deviation of $M$ from being a Cohen-Macaulay module. Here is a
summary about results we use in the sequel.

\begin{proposition} \label{2.y} Let $M$ denote a $d$-dimensional $A$-module.
Let $k \in \mathbb N$ an integer. Then the following results are
true:
\begin{itemize}
\item[(a)] $\dim K^i(M) \leq i$ for all $0 \leq i < d$ and $\dim
K(M) = d.$

\item[(b)] $\Ass K(M) = (\Ass M)_d.$

\item[(c)] $(\Ass K^i(M))_i = (\Ass M)_i$ for all $0 \leq i < d.$

\item[(d)] $K(M)$ satisfies $S_2.$

\item[(e)] $M$ satisfies $S_k$ if and only if $\dim K^i(M) \leq i
-k$ for all $0 \leq i < \dim M.$
\end{itemize}
\end{proposition}

For a finitely generated $R$-module $X$ let $(\Ass X)_i = \{
\mathfrak p \in \Ass X : \dim R/\mathfrak p = i \}$ for an integer
$i \in \mathbb Z.$ Cf. \cite[Section 1]{pS1} for the details of
the proof of Proposition \ref{2.y}.

\subsection{On commutative algebra}
Let $M$ be a finitely generated $R$-module, $R$ a commutative
Noetherian ring. Let $\Ass_R M = \{\mathfrak p_1, \ldots,
\mathfrak p_t\}$ denote the set of associated prime ideals. Let
\[
0 = Z(\mathfrak p_1) \cap \ldots \cap Z(\mathfrak p_t)
\]
denote a minimal primary decomposition of $M.$ That is,
$M/Z(\mathfrak p_i), i = 1,\ldots, t,$ is a non-zero $\mathfrak
p_i$-coprimary $R$-module.

Next we want to prove a constructive version of a result of
N.~Bourbaki (cf. \cite[Ch. IV, \S 2, Prop. 6]{nB}).

\begin{lemma} \label{2.2} With the previous notation let $S =
\{\mathfrak p_1, \ldots, \mathfrak p_s\}$ denote a subset of
$\Ass_R M$ for a certain numeration of the associated prime ideals
of $M.$ Put $U = \cap_{i=1}^s Z(\mathfrak p_i).$ Then
\[
\Ass_R M/U = S \quad {\rm and } \quad \Ass_R U = \Ass_R M
\setminus S.
\]
\end{lemma}

\begin{proof} Let $\Ass_R M = \{\mathfrak p_1, \ldots, \mathfrak p_t\}$
and $0 = Z(\mathfrak p_1) \cap \ldots \cap Z(\mathfrak p_t)$ a
minimal primary decomposition. First it is clear that $\Ass_R M/U
= S.$ Remember that $U = \cap_{i=1}^s Z(\mathfrak p_i)$ is a
reduced minimal primary decomposition. Define $V = \cap_{i=s+1}^t
Z(\mathfrak p_i).$ In order to show the second part of the claim
it will be enough to prove that $\Ass_R U = \{\mathfrak p_{s+1},
\ldots, \mathfrak p_t\}.$

First note that $ U \simeq U+V/V \subseteq M/V.$ Therefore $\Ass_R
U \subseteq \{\mathfrak p_{s+1}, \ldots, \mathfrak p_t\}$ as
easily seen. Now let $\mathfrak p \in \{\mathfrak p_{s+1}, \ldots,
\mathfrak p_t\}$ be a given prime ideal. Then $U/U \cap
Z(\mathfrak p) \simeq U+Z(\mathfrak p)/Z(\mathfrak p)$ is a
non-zero $\mathfrak p$-coprimary module. Since $U \cap Z(\mathfrak
p)$ is part of a minimal reduced primary decomposition of $0$ in
$U$ it follows that $\mathfrak p \in \Ass_R U,$ as required.
\end{proof}

\section{On the definition of formal cohomology}
\subsection{The basic definitions}
Let $(R, \mathfrak m, k)$ be a local Noetherian ring. Let $\xx =
x_1,\ldots,x_r$ denote a system of elements of $R$ and $ \mathfrak
b = \Rad(\xx R).$ Let $\Cech$ denote the \v{C}ech
complex of $R$ with respect to $\xx.$ For an $R$-module $M$ and an
ideal $\mathfrak a$ the projective system of $R$-modules
$\{M/\mathfrak a^nM\}_{n \in \mathbb N}$ induces a projective
system of $R$-complexes $\{\Cech \otimes M/\mathfrak a^nM\}.$ Its
projective limit $\varprojlim (\Cech \otimes M/\mathfrak a^nM)$ is
the main object of our investigations.

\begin{definition} \label{3.0} For an integer $i \in \mathbb Z$
the cohomology module $H^i(\varprojlim (\Cech \otimes M/\mathfrak
a^nM))$ is called the $i$-th $\mathfrak a$-formal cohomology with
respect to $\mathfrak b.$ In the case of $\mathfrak b = \mathfrak
m$ we speak simply about the $i$-th $\mathfrak a$-formal
cohomology. By abuse of notation we say also formal cohomology in
case there will be no doubt on $\mathfrak a.$
\end{definition}

In the following let $\lam = \varprojlim (\cdot \otimes
R/\mathfrak a^n)$ denote the $\mathfrak a$-adic completion. For an
$R$-module $M$
it turns out that the complex $\varprojlim (\Cech
\otimes M \otimes R/\mathfrak a^n )$ is isomorphic to $\lam (\Cech \otimes M).$ In the
derived category this complex is isomorphic to $\lam ({\rm
\Gamma}_{\mathfrak b}(I)),$ where $M \qism I$ denotes an injective
resolution of $M.$ For further results in this direction see
\cite{pS3}.

As a first result here there is a relation of the formal
cohomology with respect to the projective limits of certain local
cohomology modules.

\begin{proposition} \label{3.1} With the previous notation there
is the following short exact sequence
\[
0 \to \varprojlim{}^1 H^{i+1}_{\mathfrak b}(M/\mathfrak a^nM) \to
H^i(\varprojlim (\Cech \otimes M/\mathfrak a^nM)) \to \varprojlim
H^i_{\mathfrak b}(M/\mathfrak a^nM) \to 0
\]
for all $i \in \mathbb Z.$ In the case of $\mathfrak b = \mathfrak m$
and a finitely generated $R$-module $M$ it provides isomorphisms
\[
H^i(\varprojlim (\Cech \otimes M/\mathfrak a^nM)) \simeq
\varprojlim H^i_{\mathfrak m}(M/\mathfrak a^nM)
\]
for all $i \in \mathbb Z.$
\end{proposition}

\begin{proof} The \v{C}ech complex $\Cech$ is a complex of flat
$R$-modules. Whence the natural epimorphism $M/\mathfrak a^{n+1}M
\to M/\mathfrak a^nM, n \in \mathbb N,$ induces an $R$-morphism of
$R$-complexes
\[
\Cech \otimes M/\mathfrak a^{n+1}M \to \Cech \otimes M/\mathfrak
a^nM
\]
which is degree-wise an epimorphism. By the definition of the
projective limit there is a short exact sequence of complexes
\[
0 \to \varprojlim (\Cech \otimes M/\mathfrak a^nM) \to \prod
(\Cech \otimes M/\mathfrak a^nM) \to \prod (\Cech \otimes
M/\mathfrak a^nM) \to 0
\]
(cf. e.g. \cite{pS3}). Now the long exact cohomology sequence
provides the first part of the claim. To this end break it up into
short exact sequences and take into account that homology commutes
with direct products.

For the proof of the second part remember that $\HH^i(M/\mathfrak
a^nM), i \in \mathbb Z,$ is an Artinian $R$-module whenever $M$ is
a finitely generated $R$  (cf. \cite[Section 6]{aG}). So the
corresponding projective system satisfies the Mittag-Leffler
condition. That is, $\varprojlim{}^1$ vanishes on the projective
system of Artinian $R$-modules. The proof is now a consequence of
the first part.
\end{proof}

Let $(\hat R, \hat{\mathfrak m})$ denote the $\mathfrak m$-adic
completion of $(R, \mathfrak m).$ An Artinian $R$-module $A$ has a
natural structure of an $\hat R$-module such that the
natural homomorphisms $A \to \hat A$ and $A \to A \otimes \hat R$
are isomorphisms.

\begin{proposition} \label{3.3} Let $M$ be a finitely generated $R$-module.
Then $\varprojlim H^i_{\mathfrak m}(M/\mathfrak a^nM), i \in
\mathbb Z,$ has a natural structure as an $\hat{R}$-module and and
there are isomorphisms
\[
\varprojlim \HH^i(M/\mathfrak a^nM) \simeq \varprojlim H^i_{\hat{
\mathfrak m}}(\hat{ M}/\mathfrak a^n \hat{M})
\]
for all $i \in \mathbb Z.$
\end{proposition}

\begin{proof} Let $N$ be a finitely generated $R$-module. Then it is
known that $\HH^i(N), i \in \mathbb Z,$ is an Artinian $R$-module
(cf. e.g. \cite[Section 6]{aG}). Because of the previous remarks
and the flatness of $\hat{R}$ over $R$ there are $R$-isomorphisms
$\HH^i(N) \simeq H^i_{\hat{\mathfrak m}}(\hat {N})$ for all $i \in
\mathbb Z.$ Now take $N = M/\mathfrak a^nM$ and pass to the
projective limit. Then this proves the claim.
\end{proof}

The previous result has the advantage that one might assume the
existence of a dualizing complex in order to consider the formal
cohomology. Note that by the Cohen Structure theorem $\hat R$ is
the factor ring of a regular local ring.

Let $U = \Spec R \setminus \{\mathfrak m\}.$ Let $({\hat U},
\mathcal O_{\hat U})$ denote the formal completion of $U$ along
$V(\mathfrak a) \setminus \{\mathfrak m\}$ (cf. \cite{gF} and
\cite{PS} for the details). For an $R$-module $M$ let $\mathcal F$
denote the associated sheaf on $U.$ Let $\hat {\mathcal F}$ denote
the coherent ${\mathcal O}_{\hat U}$-sheaf associated to
$\varprojlim M/\mathfrak a^ nM.$  Let ${\hat M}^{\mathfrak a}$
denote the $\mathfrak a$-adic completion of $M.$ Moreover
$\mathcal J$ denotes the ideal sheaf of $\mathfrak a$ on $(U,
\mathcal O_{U}).$ Then there is the following relation to the
formal local cohomology (cf. also \cite{PS}).

\begin{lemma} \label{3.4} Let $M$ denote a finitely generated $R$-module.
With the previous notation there are an exact sequence
\[
0 \to \varprojlim \HH^0(M/\mathfrak a^nM) \to {\hat M}^{\mathfrak
a} \to H^0({\hat U}, \hat{\mathcal F}) \to \varprojlim
\HH^1(M/\mathfrak a^nM) \to 0
\]
and isomorphisms
\[
H^i({\hat U},\hat{\mathcal F} ) \simeq \varprojlim
\HH^{i+1}(M/\mathfrak a^nM)
\]
for all $i \geq 1.$
\end{lemma}

\begin{proof} Let $n \in \mathbb N$ denote an integer. First
remember that there is a functorial exact sequence
\[
0 \to \HH^0(M/\mathfrak a^ nM) \to M/\mathfrak a^ nM \stackrel{\phi_n}{\longrightarrow} H^0(U,
\mathcal F/\mathcal J^n \mathcal F) \to \HH^1(M/\mathfrak a^nM) \to 0
\]
and isomorphisms $H^i(U,\mathcal F/\mathcal J^n \mathcal F) \simeq
\HH^{i+1}(M/\mathfrak a^ nM)$ for all $i \in \mathbb Z$ (cf. e.g.
\cite{aG}). The family of $R$-modules $\{\im \phi_n\}_{n \in
\mathbb N},$ as a surjective system, and the families
$\{\HH^i(M/\mathfrak a^ nM)\}_{n \in \mathbb N}, i \in \mathbb N,$
as families of Artinian $R$-modules, both satisfy the
Mittag-Leffler condition. Therefore, the above exact sequence
induces -- by passing to the projective limit -- an exact sequence
\[
0 \to \varprojlim \HH^0(M/\mathfrak a^nM) \to {\hat M}^{\mathfrak
a} \to H^0({\hat U}, \hat{\mathcal F}) \to \varprojlim
\HH^1(M/\mathfrak a^nM) \to 0,
\]
which proves the first part of the claim.

The above isomorphisms provide an isomorphism
\[
\varprojlim H^i(U,\mathcal F/\mathcal J^n \mathcal F) \simeq
\varprojlim \HH^{i+1}(M/\mathfrak a^ nM)
\]
for all $i \in \mathbb Z.$ Now the natural homomorphism $H^i({\hat
U},\hat{\mathcal F} ) \to \varprojlim H^i(U,\mathcal F/\mathcal
J^n \mathcal F), i \in \mathbb Z,$ yields an isomorphism (cf.
\cite[Ch. III, Prop.2.1]{PS}). This finishes the proof of the
statement.
\end{proof}

\subsection{On duality}
In this subsection let $(R, \mathfrak m)$ denote a local ring
possessing a dualizing complex $D^{\cdot}_R.$ The main goal of the
considerations here is an expression of the formal cohomology in
terms of a certain local cohomology of the dualizing complex. To
be more precise the following result holds.

\begin{theorem} \label{3.5} Let $M$ be a finitely generated
$R$-module. For an ideal $\mathfrak a$ of $R$ there are isomorphisms
\[
\varprojlim \HH^i(M/\mathfrak a^nM) \simeq
\Hom_R(H^{-i}_{\mathfrak a}(\Hom_R(M, D^{\cdot}_R)),E),
\]
for all $i \in \mathbb Z,$ where $E = E_R(R/\mathfrak m)$ denotes
the injective hull of the residue field $k.$
\end{theorem}

\begin{proof} Let $n \in \mathbb N$ be an integer. By virtue of the Local
Duality Theorem (cf. \ref{2.1}) there are the isomorphisms
\[
\HH^i(M/\mathfrak a^nM) \simeq \Hom_R(H^{-i}(\Hom_R(M/\mathfrak
a^nM, D^{\cdot}_R)), E)
\]
for all $i \in \mathbb Z.$ By passing to the projective limit
there are isomorphisms
\[
\varprojlim \HH^i(M/\mathfrak a^nM) \simeq
\Hom_R(H^{-i}(\varinjlim \Hom_R(M/\mathfrak a^nM, D^{\cdot}_R)),
E)
\]
for all $i \in \mathbb Z.$ To this end remember that the injective
limit commutes with cohomology and is transformed into a
corresponding projective system by $\Hom$ in the first place. Now
the proof turns out because $\varinjlim \Hom(M/\mathfrak a^nM,
D^{\cdot}_R) \simeq \gam (\Hom_R(M,D^{\cdot}_R))$ as easily seen.
\end{proof}

\begin{remark} \label{2.4} In the case the local ring $(R,
\mathfrak m)$ possesses a dualizing complex it is a quotient of a
local Gorenstein ring $(S, \mathfrak n)$ (cf. \cite{tK}).
Therefore, we may use
\[
D_R^{\cdot} = \Hom_S(R, I_S^{\cdot}) [-n], \; n = \dim S,
\]
as the (normalized) dualizing complex, where $I_S^{\cdot}$ denotes
the minimal injective resolution of $S$ as an $S$-module. By the
local duality (cf. \ref{2.1})
\[
\varprojlim \HH^i(M/\mathfrak a^nM) \simeq \Hom_R(\varinjlim
\Ext_S^{n-i}(M/\mathfrak a^nM,S),E)
\]
for all $i \in \mathbb Z,$ where $E$ denotes the injective hull of
the residue field. In his unpublished habilitation (cf. \cite{jH})
Herzog introduced
\[
H^i_{\mathfrak a}(M,N) = \varinjlim \Ext_R^i(M/\mathfrak a^nM, N),
i \in \mathbb Z,
\]
for two $R$-moduls $M,N$ and an ideal $\mathfrak a \subset R$ as
the {\sl generalized local cohomology} with respect to $\mathfrak
a.$ With the previous notation there are isomorphisms
\[
\varprojlim \HH^i(M/\mathfrak a^nM) \simeq \Hom_R(H^{n-i}_{\mathfrak a
S}(M,S), E), \; i \in \mathbb Z,
\]
where $M$ is considered as an $S$-module. So, the $i$-th
$\mathfrak a$-formal cohomology $\varprojlim \HH^i(M/\mathfrak
a^nM)$ is isomorphic to the Matlis dual of $H^{n-i}_{\mathfrak a
S}(M,S)$ equipped with its natural $R$-module structure.
\end{remark}

The previous result has as a consequence a non-vanishing behavior
of the formal cohomology, important for the subsequent
considerations.

\begin{corollary} \label{3.6} Let $\mathfrak p$ denote a prime
ideal and $i \in \mathbb Z$ be such that $\varprojlim
H^i_{\mathfrak p R_{\mathfrak p}}(M_{\mathfrak p}/\mathfrak a^n
M_{\mathfrak p}) \not= 0.$ Then $\varprojlim \HH^{i+\dim
R/\mathfrak p}(M/\mathfrak a^nM) \not= 0.$
\end{corollary}

\begin{proof} By virtue of Matlis' duality for the local ring
$R_{\mathfrak p}$ it follows that $H^{-i}_{\mathfrak a
R_{\mathfrak p}}(\Hom (M_{\mathfrak p}, D^{\cdot}_{R_{\mathfrak
p}}))$ does not vanish (cf. \ref{3.5}). Now there is an
isomorphism of complexes
\[
\Hom(M_{\mathfrak p}, D^{\cdot}_{R_{\mathfrak p}}) \simeq
\Hom_R(M, D^{\cdot}_R)[-\dim R/\mathfrak p] \otimes R_{\mathfrak
p}
\]
(cf. \ref{2.1} and remember that $M$ is a finitely generated
$R$-module). But this provides the isomorphisms
\[
H^{-j}_{\mathfrak a R_{\mathfrak p}}(\Hom (M_{\mathfrak p},
D^{\cdot}_{R_{\mathfrak p}})) \simeq H^{-j-\dim R/\mathfrak
p}_{\mathfrak a}(\Hom (M, D^{\cdot}_R)) \otimes R_{\mathfrak p}
\]
for all $j \in \mathbb Z.$ Therefore $H^{-i-\dim R/\mathfrak
p}_{\mathfrak a}(\Hom_R(M, D^{\cdot}_R)) \not= 0.$ By Matlis'
duality this implies the non-vanishing of $\varprojlim \HH^{i+\dim
R/\mathfrak p}(M/\mathfrak a^nM)$  (cf. \ref{3.5}). This completes
the proof.
\end{proof}

We conclude this subsection with the proof of the fact that
equivalent ideal topologies define isomorphic formal cohomology
modules. Here $\{M_n\}_{n \in \mathbb N}$ is called a decreasing
family of submodules provided $M_{n+1} \subseteq M_n$ for all $n
\in \mathbb N.$

\begin{lemma} \label{3.6a} Let $M$ be a finitely generated $R$-module.
Let $\{M_n\}_{n \in \mathbb N}$ be a decreasing family
of submodules of $M.$  Suppose that their topology is
equivalent to the $\mathfrak a$-adic topology on $M.$
Then there are isomorphisms
\[
\varprojlim \HH^i(M/\mathfrak
a^nM) \simeq \varprojlim \HH^i(M/M_n)
\]
for all $i \in \mathbb Z.$
\end{lemma}

\begin{proof} Let $\Cech$ denote the \v{C}ech complex of $R$ with respect
to a system of elements $\xx$ such that $\Rad \xx R = \mathfrak
m.$ For any flat $R$-module $F$ there is an isomorphism
$\varprojlim F \otimes (M/M_n) \simeq \varprojlim F \otimes
M/\mathfrak a^ nM.$ To this end remember that $F \otimes (M/N)
\simeq (F \otimes M)/(F \otimes N)$ for any submodule $N \subseteq
M.$ Moreover, $\{F \otimes M_n\}$ is equivalent to the $\mathfrak
a$-adic topology on $F\otimes M.$

Since $\Cech$ is a bounded complex of flat $R$-modules this isomorphism  extends to an
isomorphism $\varprojlim \Cech \otimes M/M_n \simeq \varprojlim \Cech \otimes
M/\mathfrak a^n M $ of $R$-complexes. Therefore, it will be enough to show that
\[
H^i(\varprojlim \Cech \otimes (M/M_n)) \simeq \varprojlim H^i(\Cech \otimes M/M_n),
\; i \in \; \mathbb Z,
\]
(cf. \ref{3.1}). Since $H^i(\Cech \otimes (M/M_n)) \simeq \HH^i(M/M_n), i \in \mathbb Z,$
is an Artinian $R$-module this follows by the Mittag-Leffler arguments
as in the proof of the second part of \ref{3.1}.
\end{proof}

As a first structure result on the formal cohomology modules
$\varprojlim \HH^i(M/\mathfrak a^nM), i \in \mathbb Z,$ for a
finitely generated $R$-module $M$ we consider their behavior with
respect to the $\mathfrak a$-adic completion. Let ${\rm L}_i \lam,
i \in \mathbb Z,$ denote the left derived functors of the
$\mathfrak a$-adic completion functor $\varprojlim (\cdot \otimes
R/\mathfrak a^n)$ (cf. \cite{GM}, \cite{aS} for the basic results
for modules and \cite{pS3} for an extension to complexes). An
extensive consideration of the functors ${\rm L}_i \lam, i \in
\mathbb Z,$ has been done in the fundamental work \cite{ALL}.

\begin{theorem} \label{3.7a} Let $\mathfrak a$ denote an ideal of an
arbitrary local ring $(R, \mathfrak m).$ Let $M$ be a finitely generated $R$-module.
For an integer $j \in \mathbb Z$ there are the following isomorphisms
\begin{equation*}
{\rm L}_i \lam( \varprojlim \HH^i(M/\mathfrak a^nM)) \simeq
    \begin{cases} 0 & \text{ for } i \not= 0,\\
                  (\varprojlim \HH^j(M/\mathfrak a^n M))^{\mathfrak a} & \text{ for } i = 0.
    \end{cases}
\end{equation*}
Moreover, $\varprojlim \HH^j(M/\mathfrak a^nM)$ is an $\mathfrak
a$-adic complete $R$-module, i.e. $(\varprojlim \HH^j(M/\mathfrak
a^n M))^{\mathfrak a} \simeq \varprojlim \HH^j(M/\mathfrak a^n
M)).$
\end{theorem}

\begin{proof} Without loss of generality we may assume that $(R, \mathfrak m)$ admits
a dualizing complex $D^{\cdot}_R$ (cf. \ref{3.3}). For simplicity of notation
put $X^j := \varprojlim \HH^j(M/\mathfrak a^n M), j \in \mathbb Z.$ Then there is
the following isomorphism $X^j \simeq \Hom (H^j, E),$ where $H^j := H^{-j}_{\mathfrak a}(\Hom (M, D^{\cdot}_R)),$
(cf. \ref{3.5}).

Let $X$ denote an $R$-module. For the computation of ${\rm L}_i \lam (X), i \in
\mathbb Z,$ there is the following short exact sequence
\[
0 \to \varprojlim{}^1 \Tor_{i+1}^R(R/\mathfrak a^n, X) \to {\rm L}_i \lam (X) \to
\varprojlim \Tor_i^R(R/\mathfrak a^n, X) \to 0
\]
(cf. \cite[Prop. 1.1]{GM} or \cite{pS3}). Thus, for the first part of our claim
it will be enough to prove that $\varprojlim{}^1 \Tor_{i+1}^R(R/\mathfrak a^n, X) = 0$
for all $i \in \mathbb Z$ and $\varprojlim \Tor_i^R(R/\mathfrak a^n, X) = 0$ for all
integers $i \not= 0.$

To this end consider $H^i_{\mathfrak a}(H^j) \simeq \varinjlim \Ext^i(R/\mathfrak a^n, H^j).$
Because of $\Supp H^j \subseteq V(\mathfrak a)$ clearly $H^i_{\mathfrak a}(H^j) = 0$
for all $i \not= 0$ and $H^0_{\mathfrak a}(H^j) \simeq H^j.$ By the definition of the
direct limit there is the following, canonical exact sequence
\[
0 \to \bigoplus_{n \in \mathbb N} \Ext^i(R/\mathfrak a^n, H^j) \stackrel{\Phi_i}{\longrightarrow}
\bigoplus_{n \in \mathbb N} \Ext^i(R/\mathfrak a^n, H^j) \to H^i_{\mathfrak a}(H^j) \to 0,
\]
where $\Phi_i$ is defined by the definition of the direct limit.
Now apply the Matlis duality functor $\Hom( \cdot, E).$ Because of
\[
\Hom( \Ext^i(R/\mathfrak a^n, H^j), E) \simeq \Tor_i(R/\mathfrak a^n, X^j)
\]
for all $i, j \in \mathbb Z$ and all $n \in \mathbb N,$ it
transforms the direct system $\{ \Ext^i(R/\mathfrak a^n, H^j)\}_{n
\in \mathbb N}$ into the inverse system $\{\Tor_i(R/\mathfrak a^n,
X^j)\}_{n \in \mathbb N}.$ Moreover it provides the short exact
sequences
\[
0 \to \Hom(H^i_{\mathfrak a}(H^j), E) \to \prod_{n \in \mathbb N} \Tor_i(R/\mathfrak a^n, X^j)
\stackrel{\Psi_i}{\longrightarrow} \prod_{n \in \mathbb N} \Tor_i(R/\mathfrak a^n, X^j) \to 0
\]
for all $i, j \in \mathbb Z.$ By the definition of the homomorphism $\Psi_i$ it follows that
\[
\Coker \Psi_i \simeq \varprojlim{}^1 \Tor_i^R(R/\mathfrak a^n, X) \quad \mbox{ and } \quad
\Ker \Psi_i \simeq \varprojlim \Tor_i^R(R/\mathfrak a^n, X).
\]
By the vanishing of the local cohomology of $H^j$ this provides the vanishing results of $\varprojlim{}^1$
and $\varprojlim$ of the $\Tor$'s as claimed above. Moreover, for $i = 0$ it yields the
isomorphisms
\[
X^j \simeq \Ker \Psi_0 \simeq (\varprojlim \HH^i(M/\mathfrak a^nM))^{\mathfrak a}.
\]
To this end remember that $X^j \simeq \Hom( H^j,E),$ as mentioned
above. This finally completes the proof of the result.
\end{proof}

The class ${\mathcal C}_{\mathfrak a}$ of $R$-modules $X$ such that ${\rm L}_i \lam(X) = 0$
for $i > 0$ and $  {\rm L}_0 \lam(X) = X^{\mathfrak a}$ has been introduced by Simon (cf.
\cite[5.2]{aS}). Therefore, the $\mathfrak a$-formal cohomology modules of a finitely generated
$R$-module $M$ belong to ${\mathcal C}_{\mathfrak a}.$

As a corollary there is the following Nakayama type criterion
about the vanishing of the $\mathfrak a$-formal cohomology.

\begin{corollary} \label{3.7b} Let $M$ denote a finitely generated $R$-module. Let $j \in
\mathbb Z.$ Suppose that $\varprojlim \HH^j(M/\mathfrak a^nM) =
\mathfrak a (\varprojlim \HH^j(M/\mathfrak a^nM)).$ Then
$\varprojlim \HH^j(M/\mathfrak a^nM) = 0.$
\end{corollary}

\begin{proof} For simplicity of notation put $ \varprojlim \HH^j(M/\mathfrak a^nM) = X.$
The assumption provides $X = \mathfrak a^n X, n \in \mathbb N,$ as follows by an induction.
Therefore
\[
0 = \varprojlim X/\mathfrak a^n X = X^{\mathfrak a}.
\]
By the Theorem \ref{3.7a} $X$ is $\mathfrak a$-adically complete. Therefore $X = X^{\mathfrak a}$
and $X = 0,$ as required.
\end{proof}

\subsection{Exact sequences}
First of all we want to relate the behavior of the formal
cohomology with respect to short exact sequences of $R$-modules.
This is a technical tool that simplifies arguments in further
considerations.

\begin{theorem} \label{3.6b} Let $(R, \mathfrak m)$ denote a local ring. Let
$0 \to A \to B \to C \to 0$ denote a short exact sequence of finitely generated $R$-modules.
For an ideal $\mathfrak a$ of $R$ there is a long exact sequence
\[
\ldots \to \varprojlim \HH^i(A/\mathfrak a^n A) \to \varprojlim \HH^i(B/\mathfrak a^n B)
\to \varprojlim \HH^i(C/\mathfrak a^n C) \to \varprojlim \HH^{i+1}(A/\mathfrak a^n A)  \to \ldots
\]
\end{theorem}

\begin{proof} For any finitely generated $R$-module $M$ the formal cohomology
of $M$ and $\hat{M}$ coincide (cf. \ref{3.1}). So we may assume the existence of a
dualizing complex $D^{\cdot}_R.$ Let $\Cech$ denote the \v{C}ech complex of
$R$ with respect to a system of elements $\xx$ such that $\Rad \xx R = \Rad \mathfrak a.$
The short exact sequence $0 \to A \to B \to C \to 0$ induces a short exact
sequence of $R$-complexes
\[
0 \to \Cech \otimes \Hom (C, D^{\cdot}_R) \to \Cech \otimes \Hom (B, D^{\cdot}_R)
\to \Cech \otimes \Hom (A, D^{\cdot}_R) \to 0.
\]
Remember that $D^{\cdot}_R$ resp. $\Cech$ is a bounded complex of
injective resp. flat $R$-modules. By passing to the Matlis dual
and taking the long exact cohomology sequence this proves the
claim. Remember that $H^i(\Cech \otimes \Hom(M,D^{\cdot}_R))
\simeq H^i_{\mathfrak a}(\Hom(M, D^{\cdot}_R))$ for all $i \in
\mathbb Z.$
\end{proof}

\begin{remark} \label{3.6c} One might ask for a corresponding result
for a short exact sequence $0 \to A \to B \to C \to 0,$ where the
$R$-modules are not necessarily finitely generated. It is not
clear whether this will be true.

An alternative proof of \ref{3.6b} works as follows. The short exact sequence induces
a projective system of short exact sequences
\[
0 \to \Cech \otimes A/B \cap \mathfrak a^n A \to \Cech B/\mathfrak a^n B
\to \Cech \otimes C/\mathfrak a^n C \to 0
\]
for all $n \in \mathbb N.$ Because $\Cech$ is a complex of flat
$R$-modules and because the maps
\[
A/B\cap \mathfrak a^{n+1} A \to A/B \cap \mathfrak a^n A
\]
are surjective it follows that the projective system of
$R$-complexes $\{ \Cech \otimes A/B \cap \mathfrak a^n A \}_{n \in
\mathbb Z}$ satisfies degree-wise the Mittag-Leffler condition.
Therefore the projective limit provides a short exact sequence of
complexes
\[
0 \to \varprojlim \Cech \otimes A/B \cap \mathfrak a^n A \to \varprojlim \Cech \otimes B/\mathfrak a^n B
\to \varprojlim \Cech \otimes C/\mathfrak a^n C \to 0.
\]
By view of the long exact cohomology sequence it follows (cf. the definition and \ref{3.1})
that there a long exact sequence
\[
\ldots \to \varprojlim \HH^i(A/B \cap \mathfrak a^ n A) \to \varprojlim \HH^i(B/\mathfrak a^n B) \to
\varprojlim \HH^i(C/\mathfrak a^ n C) \to \ldots .
\]
In the case $\{B \cap \mathfrak a^ n A\}$ is equivalent to the
$\mathfrak a$-adic topology on $A$ this yields another proof of
the exact sequence in \ref{3.6b} (cf. \ref{3.6a}). By the
Artin-Rees Lemma (cf. \cite[Ch. III, \S 3, Cor. 1]{nB}) this is
true in case $B$ is a finitely generated $R$-module.
\end{remark}

As an application let us consider the behavior
of the formal cohomology by factoring out the $\mathfrak
m$-torsion.

\begin{corollary} \label{3.10} Let $(R, \mathfrak m)$ denote a local
ring. For a finitely generated $R$-module $M$ let $N \subseteq M$
be an $R$-module such that $\Supp N \cap V(\mathfrak a) \subseteq
V(\mathfrak m)$. Put $\bar{M} = M/N.$ Then there is a short exact
sequence
\[
0 \to N^{\mathfrak a} \to \varprojlim \HH^0(M/\mathfrak a^n M) \to \varprojlim
\HH^0(\bar{M}/\mathfrak a^n \bar{M}) \to 0
\]
and isomorphisms $\varprojlim \HH^i(M/\mathfrak a^n M)\simeq
\varprojlim \HH^i(\bar{M}/\mathfrak a^n \bar{M})$ for all
$i \geq 1.$
\end{corollary}

\begin{proof} There is the following short exact sequence $0 \to N \to M
\to \bar{M} \to 0.$ Then there is the long exact sequence
\begin{multline*}
\ldots \to \varprojlim \HH^i(N/\mathfrak a^n N) \to \varprojlim \HH^i(M/\mathfrak a^n M)
\to \\ \to
\varprojlim \HH^i(\bar{M}/\mathfrak a^n \bar{M}) \to \varprojlim \HH^{i+1}(N/\mathfrak a^n N)
\to \ldots
\end{multline*}
(cf. \ref{3.6b}). By view of the assumption $\Supp N \cap
V(\mathfrak a) \subseteq V(\mathfrak m)$ it follows that
$N/\mathfrak a^ nN$ is an $R$-module of finite length for all $n
\in \mathbb N.$ That is, $\HH^i(N/\mathfrak a^ nN) = 0$ for $i >
0$ and all $n \in \mathbb N.$ Moreover $\HH^0(N/\mathfrak a^ nN)
\simeq N/\mathfrak a^n N$ and therefore $\varprojlim
\HH^0(N/\mathfrak a^n N) \simeq N^{\mathfrak a}.$ So the above
long exact sequence provides the short exact sequence and the
isomorphisms of the claim.
\end{proof}

In the subsequent section there is a generalization of \ref{3.10}.
In fact there is a precise computation of the $0$-th formal
cohomology.

\begin{theorem} \label{3.7} Let $M$ be a finitely generated
$R$-module. Choose $x \in \mathfrak m$ an element such that $x
\not\in \mathfrak p$ for all $\mathfrak p \in \Ass_R M \setminus
\{\mathfrak m\}.$ Then there are short exact sequences
\[
0 \to H_0(x; \varprojlim \HH^i(M/\mathfrak a^nM)) \to \varprojlim
H^i_{\mathfrak m}(M'/\mathfrak a^nM') \to H_1(x; \varprojlim
\HH^{i+1}(M/\mathfrak a^nM)) \to 0
\]
for all $i \in \mathbb Z,$ where $M' = M/xM.$
\end{theorem}

\begin{proof} By the choice of $x$ it follows that $0 :_M x$
is an $R$-module of finite length. Moreover the multiplication by $x$ induces
an exact sequence
\[
0 \to 0:_M x \to M \stackrel{x}{\longrightarrow} M \to M' \to 0
\]
breaks into two short exact sequences $0 \to N \to M \to \bar{M}
\to 0,$ where $N = 0 :_M x$ and $\bar{M} = M/N,$ and $0 \to
\bar{M} \stackrel{x}{\to} M \to M' \to 0.$

The first of these sequences induces isomorphisms $\varprojlim \HH^i(M/\mathfrak a^nM)
\simeq \varprojlim \HH^i(\bar{M}/\mathfrak a^n \bar{M})$ for all $i > 0$ and a short exact sequence
\[
0 \to N \to \varprojlim \HH^0(M/\mathfrak a^n M) \to \varprojlim
\HH^0(\bar{M}/\mathfrak a^n \bar{M}) \to 0
\]
(cf \ref{3.10}). The second sequence induces a long exact sequence for the formal cohomology
modules
\[
\ldots \to \varprojlim \HH^i(\bar{M}/\mathfrak a^n \bar{M})
\stackrel{x}{\to} \varprojlim \HH^i(M/\mathfrak a^n M) \to
\varprojlim \HH^i(M'/\mathfrak a^n M')  \to \ldots
\]
(cf \ref{3.6b}).  With the isomorphisms above this proves the claim
for $i > 0.$ To this end one has to break up the long exact sequence
into short exact sequences.

For the proof in the case $i = 0,$ the only remaining case, consider the composite of
the above short exact sequence with the previous one for $i= 0.$ Then this
completes the proof for $i = 0.$
\end{proof}

Another short exact sequence relates the $\mathfrak a$-formal
cohomology to the $(\mathfrak a, xR)$-formal cohomology for any
element $x \in \mathfrak m.$ To be more precise:

\begin{theorem} \label{3.8} Let $x \in \mathfrak m$ denote an
element of $(R, \mathfrak m).$ For an ideal $\mathfrak a$ and a
finitely generated $R$-module $M$ there is the long exact sequence
\[
\ldots \to \Hom(R_x,\varprojlim \HH^i(M/\mathfrak a^nM)) \to
\varprojlim \HH^i(M/\mathfrak a^nM) \to \varprojlim
\HH^i(M/(\mathfrak a,x)^nM) \to \ldots
\]
for all $i \in \mathbb Z.$
\end{theorem}

\begin{proof} Without loss of generality (cf. \ref{3.3}) we may
assume that $R$ admits a dualizing complex $D^{\cdot}_R.$ The
\v{C}ech complex $\check{C}_x$ of the single element $x$ is the
fibre of the natural homomorphism $R \to R_x.$ So there is a split
exact sequence
\[
0 \to R_x[-1] \to \check{C}_x \to R \to 0.
\]
Let $\xx$ denote a system of elements of $R$ such that $\Rad
\mathfrak a = \Rad \xx R.$ By tensoring the above short exact
sequence of flat $R$-modules with $\Cech \otimes \Hom(M,
D^{\cdot}_R)$ it provides an exact sequence of $R$-complexes
\[
0 \to \Cech \otimes \Hom(M, D^{\cdot}_R) \otimes R_x[-1] \to
\check{C}_{\xx ,x} \otimes \Hom(M, D^{\cdot}_R) \to \Cech \otimes
\Hom(M, D^{\cdot}_R) \to 0.
\]
Notice that the above short exact sequence of complexes is split
exact. Taking the long exact cohomology sequence it provides an
exact sequence
\[
\cdots \to H^j_{(\mathfrak a,xR)}(\Hom(M,D^{\cdot}_R)) \to
H^j_{\mathfrak a}(\Hom(M,D^{\cdot}_R))\to H^j_{\mathfrak
a}(\Hom(M,D^{\cdot}_R))\otimes R_x \to \cdots
\]
for all $j \in \mathbb Z.$ By applying Matlis' duality it provides
the exact sequence of the statement (cf. \ref{3.5}).
\end{proof}

As an application of Theorem \ref{3.8} there is an exact sequence
for the formal cohomology with respect to an ideal generated by a
single element.

\begin{corollary} \label{3.9} Let $x \in \mathfrak m$ denote an
element. Let $M$ be a finitely generated $R$-module. Then there is
a short exact sequence
\[
\ldots \to \Hom(R_x,\HH^i(M)) \to H^i_{\mathfrak m}(M) \to
\varprojlim \HH^i(M/x^nM) \to \ldots
\]
for all $i \in \mathbb Z.$
\end{corollary}

\begin{proof} The corollary is a consequence of Theorem \ref{3.8}
with the particular case $\mathfrak a = 0.$
\end{proof}

\section{Vanishing results}
\subsection{On the 0-th formal cohomology} Let $(R, \mathfrak m)$
denote a local ring. Let $M$ be a finitely generated $R$-module.
For an $R$-submodule $N$ of $M$ denote by $N :_M \langle \mathfrak
m \rangle$ the ultimate constant $R$-module $N :_M \mathfrak m^n,
n $ large.

Let $0 = \cap_{\mathfrak p \in \Ass M} Z(\mathfrak p)$ denote a
minimal primary decomposition of $0$ in $M.$ Moreover, let
$\mathfrak a$ denote an ideal of $R.$ Then define
\[
T_{\mathfrak a}(M) = \{ \mathfrak p \in \Ass_R M : \dim
R/(\mathfrak a, \mathfrak p) = 0 \}.
\]
Furthermore, put
\[
u_M(\mathfrak a) = \bigcap_{\mathfrak p \in \Ass_R M \setminus
T_{\mathfrak a}(M)} Z(\mathfrak p).
\]
Now it will be shown that $u_M(\mathfrak a)$ plays an important
r\^{o}le in order to understand the 0-th formal cohomology module.
To this end denote by $\hat R$ the completion of $R$ and $\hat{M}
\simeq M \otimes \hat R$ the completion of the finitely generated
$R$-module $M.$

\begin{lemma} \label{4.1} With the previous notation we have:
\begin{itemize}
\item[(a)] $ \bigcap_{n \geq 1} (\mathfrak a^n M :_M \langle
\mathfrak m \rangle) = u_M(\mathfrak a) .$

\item[(b)] $\Ass_R (u_M(\mathfrak a)) = T_{\mathfrak a}(M).$

\item[(c)] $\varprojlim \HH^0(M/\mathfrak a^nM) \simeq u_{\hat
M}(\mathfrak a \hat{R}).$
\end{itemize}
\end{lemma}

\begin{proof} The proof of (a) is easily seen because of
\[
\bigcap_{n \geq 1} (\mathfrak a^n M :_M \langle \mathfrak m
\rangle) = \bigcap_{\mathfrak P \in \Supp M/\mathfrak aM \setminus
V(\mathfrak m)} \ker(M \to M_{\mathfrak P})
\]
(cf. \cite[(2.1)]{pS0} for the details). Then the statement in (b)
is a consequence of (a) (cf. \ref{2.2}).

In order to proof (c) first note that one may assume $M = \hat{M}$
and $R = \hat R$ as follows by passing to the completion (cf.
\ref{3.3}). But now $\HH^0(M/\mathfrak a^nM) \simeq \mathfrak a^n
M:_M \langle \mathfrak m \rangle /\mathfrak a^nM.$ So there is a
short exact sequence of inverse systems
\[
0 \to \{\mathfrak a^nM\}_{n \in \mathbb N} \to \{\mathfrak a^n
M:_M \langle \mathfrak m \rangle\}_{n \in \mathbb N} \to
\{\HH^0(M/\mathfrak a^nM)\}_{n \in \mathbb N} \to 0.
\]
By passing to the projective limit it provides an injection
\[
0 \to \bigcap_{n \geq 1} (\mathfrak a^n M :_M \langle \mathfrak m
\rangle) \stackrel{\phi}{\to} \varprojlim \HH^0(M/\mathfrak a^nM).
\]
In order to finish it will be enough to prove that $\phi$ is
surjective. To this end let
\[\{y_n + \mathfrak a^nM\} \in
\varprojlim \HH^0(M/\mathfrak a^nM),
\]
where $y_n \in\mathfrak a^nM :_M \langle \mathfrak m \rangle$ for
all $n \in \mathbb N.$ This sequence defines an element $z \in
\varprojlim M/\mathfrak a^nM = M.$ Note that $M$ as an $\mathfrak
m$-adically complete module is also $\mathfrak a$-adically
complete (cf. \cite[Ch. VIII]{ZS}). That is, for every $n \in
\mathbb N$ there exists an $n_0 \geq n$ such that $z - y_m \in
\mathfrak a^nM$ for all $m \geq n_0.$ Therefore $z \in \cap_{m
\geq 1} (\mathfrak a^m M :_M \langle \mathfrak m \rangle),$ as
required.
\end{proof}

By view of \ref{4.1} there is the following vanishing result for
the 0-th formal cohomology.

\begin{corollary} \label{4.2} With the previous notation we have
that $\varprojlim \HH^0(M/\mathfrak a^nM) = 0$ if and only if
$\dim \hat{R}/(\mathfrak a \hat{R}, \mathfrak p) > 0$ for all
$\mathfrak p \in \Ass_{\hat{R}} \hat{M}.$

In particular, the vanishing $\varprojlim \HH^0(M/\mathfrak a^nM)
= 0$ implies that $\depth M > 0.$
\end{corollary}

\begin{proof} It turns out that $\varprojlim \HH^0(M/\mathfrak a^nM) =
0$ if and only if $\Ass_{\hat R} (u_{\hat M}(\mathfrak a \hat{R}))
= \emptyset.$ But this is equivalent to the statement (cf.
\ref{4.1}). In particular, $\varprojlim \HH^0(M/\mathfrak a^nM) =
0$ implies that $\hat{\mathfrak m} \not\in \Ass_{\hat R} \hat{M},$
whence $\depth M > 0.$
\end{proof}

Next we want to extent the statement in \ref{3.10}.

\begin{corollary} \label{4.3} Let $(R,\mathfrak m)$ denote a complete
local ring. For a finitely generated $R$-module $M$ put $U =
u_M(\mathfrak a)$ and $\bar{M} = M/U.$ Then:
\begin{itemize}
\item[(a)]$ \varprojlim \HH^0(M/\mathfrak a^n M) \simeq U$ and \,
$\varprojlim \HH^0(\bar{M}/\mathfrak a^n \bar{M}) = 0.$

\item[(b)] $ \varprojlim \HH^i(M/\mathfrak a^n M) \simeq
\varprojlim \HH^i(\bar{M}/\mathfrak a^n \bar{M}) $ for all $i
\geq1.$
\end{itemize}
\end{corollary}

\begin{proof} For the proofs of the statements in (a) see
\ref{4.1}. Now observe that
\[
\Supp U \cap V(\mathfrak a) = (\cup_{\mathfrak p \in T_{\mathfrak
a}(M)} V(\mathfrak p)) \cap V(\mathfrak a) \subseteq V(\mathfrak
m).
\]
By virtue of \ref{3.10} this proves the isomorphisms in (b).
\end{proof}

\subsection{A non-vanishing result}
The aim of this subsection will be to determine the integer
\[
\sup \{i \in \mathbb Z : \varprojlim \HH^i(M/\mathfrak a^nM) \not=
0\}.
\]
Here let $M$ be  a finitely generated $R$-module. Let $\mathfrak
a$ denote an ideal in the local ring $(R, \mathfrak m).$ We start
with an almost trivial observation.

\begin{proposition} \label{4.4} Let $\mathfrak a$ be an
ideal such that $\dim M/\mathfrak aM = 0.$ Then
\[
\varprojlim \HH^i(M/\mathfrak a^nM) \simeq \left\{
    \begin{array}{cl} 0 & \mbox{ for } i \not= 0  \mbox{ and }\\
        M^{\mathfrak a} & \mbox{ for } i = 0.
    \end{array} \right.
\]
\end{proposition}

\begin{proof} It follows that $\HH^i(M/\mathfrak a^nM) = 0$ for
all $i \not= 0.$ Notice that $M/\mathfrak aM$ is an $R$-module of
finite length. Furthermore, it provides $\HH^0(M/\mathfrak a^nM)
\simeq M/\mathfrak a^nM.$ Passing to the projective limit finishes
the proof.
\end{proof}

Now the preparation for the first non-vanishing result is
finished.

\begin{theorem} \label{4.5} Let $\mathfrak a$ denote an ideal of
$(R, \mathfrak m).$ Then
\[
\dim_R M/\mathfrak aM = \sup \{i \in \mathbb Z : \varprojlim
\HH^i(M/\mathfrak a^nM) \not= 0\}
\]
for a finitely generated $R$-module $M.$
\end{theorem}

\begin{proof} Because of $\dim M/\mathfrak a^nM = \dim M/\mathfrak a M$
for all $n \in \mathbb N$ we first note that $\HH^i(M/\mathfrak
a^nM)$ vanishes for all $i > \dim_R M/\mathfrak aM$ (cf. e.g.
\cite[Proposition 1.12]{aG}). Therefore
\[
\dim_R M/\mathfrak aM \geq \sup \{i \in \mathbb Z : \varprojlim
\HH^i(M/\mathfrak a^nM) \not= 0\}.
\]
Second note that we may assume the existence of a dualizing
complex (cf. \ref{3.3}).

In order to prove the equality take $\mathfrak p \in \Supp_R M
\cap V(\mathfrak a)$ such that $\dim R/\mathfrak p = \dim_R
M/\mathfrak aM.$ Then $\varprojlim H^0_{\mathfrak p A_{\mathfrak
p}}(M_{\mathfrak p}/\mathfrak a^n M_{\mathfrak p}) \not= 0$ (cf.
\ref{4.4}). Observe that $M_{\mathfrak p}/\mathfrak a M_{\mathfrak
p}$ is a zero-dimensional $R_{\mathfrak p}$-module. Therefore
$\varprojlim \HH^{\dim R/\mathfrak p} (M/\mathfrak a^nM) \not= 0$
(cf. \ref{3.6}).
\end{proof}

\begin{remark} \label{R x} Another proof for the non-vanishing
of $\varprojlim \HH^d(M/\mathfrak a^nM), d = \dim M/\mathfrak aM,$ can
be seen as follows. First note $\dim M/\mathfrak a^n M = d$ for all
$n \in \mathbb N.$ Then the short exact sequence
\[
0 \to \mathfrak a^nM/\mathfrak a^{n+1}M \to M/\mathfrak a^{n+1}M \to
M/\mathfrak a^nM \to 0
\]
induces an epimorphism $\HH^d(M/\mathfrak a^{n+1}M) \to \HH^d(M/\mathfrak a^nM)
\to 0,$ of non-zero $R$-modules for all $n \in \mathbb N.$ Remember that
$\dim \mathfrak a^nM/\mathfrak a^{n+1}M \leq d$ and therefore
$\HH^{d+1}(\mathfrak a^nM/\mathfrak a^{n+1}M ) = 0.$
Whence the inverse limit $\varprojlim \HH^d(M/\mathfrak a^nM)$ is not
zero.
\end{remark}

\subsection{The formal grade}
Let $M$ denote a finitely generated $R$-module, where $(R,
\mathfrak m)$ is a local ring. For an ideal $\mathfrak a$ it is
shown that $\sup \{i \in \mathbb Z : \varprojlim \HH^i(M/\mathfrak
a^nM) \not= 0\}$ is equal to $\dim_RM/\mathfrak aM$ (cf. \ref{4.5}). Now
we start to investigate the infimum for the non-vanishing.

\begin{definition} \label{4.6} For an ideal $\mathfrak a$ of $R$
define the formal grade, $\fgrade (\mathfrak a, M),$ by
\[
\fgrade (\mathfrak a, M) = \inf \{i \in \mathbb Z : \varprojlim
\HH^i(M/\mathfrak a^nM) \not= 0\}.
\]
Note that the ordinary grade is defined by $\grade (\mathfrak a, M) = \inf \{i \in
\mathbb Z : H^i_{\mathfrak a}(M) \not= 0\}$ (cf. \cite{aG}).
\end{definition}

The notion of formal grade was introduced by Peskine and Szpiro
(cf. \cite{PS}). Not so much is known about it. We continue here
with a few more investigation on the formal grade. In the
following lemma (cf. \ref{4.7}) there is a summary of basic
results.

\begin{lemma} \label{4.7} Let $\mathfrak a$ denote an ideal of
$(R, \mathfrak m).$ Let $M$ be a finitely generated $R$-module.
\begin{itemize}
\item[(a)] $\fgrade( \mathfrak a, {M/xM})  \geq \fgrade (\mathfrak a, M)
- 1,$ provided $x \not\in \mathfrak p$ for all $\mathfrak p \in
\Ass_R M \setminus \{\mathfrak m\}$.

\item[(b)] $\fgrade (\mathfrak a, M) \leq \min \{\depth_R M, \dim
M/\mathfrak aM\}.$

\item[(c)] Suppose that $R$ possesses a dualizing complex. Then
\[
\fgrade (\mathfrak a, M) \leq \fgrade ( \mathfrak a
R_{\mathfrak p}, {M_{\mathfrak p}}) + \dim R/\mathfrak p
\]
for all $\mathfrak p \in \Supp M \cap V(\mathfrak a).$

\item[(d)] Suppose that $R$ is a Gorenstein ring. Then $\fgrade (
\mathfrak a,R) + \cd ( \mathfrak a,R) = \dim R.$
\end{itemize}
\end{lemma}

\begin{proof} By virtue of the short exact sequences in \ref{3.7}
it follows that $\fgrade ( \mathfrak a, M/xM)\geq \fgrade(\mathfrak a, M)
-1.$ That is, the statement (a) is shown.

In order to prove (b) first note $\fgrade (\mathfrak a, M) \leq
\dim M/\mathfrak a M$ (cf. \ref{4.5}). Next we prove $\fgrade (
\mathfrak a, M) \leq \depth_R M$ by an induction on $t = \fgrade (
\mathfrak a, M).$ In case $t = 0$ the claim holds trivially. So
let $t \geq 1.$ Then $\varprojlim \HH^0(M/\mathfrak a^nM) = 0$ by
the definition of the formal grade. Therefore there is an
$M$-regular element $x \in \mathfrak m$ (cf. \ref{4.2}). Whence
\[
t -1 \leq \fgrade (\mathfrak a,M/xM)  \leq \depth M/xM = \depth M-1
\]
by the aid of (a) and the induction hypothesis. So the proof of
(b) is complete.

For the proof of (c) let $t = \fgrade (\mathfrak
aR_{\mathfrak p}, M_{\mathfrak p}) .$ Then $\fgrade (\mathfrak a, M) \leq t + \dim
R/\mathfrak p,$ (cf. \ref{3.6}).

Let $R$ be a Gorenstein ring. Then $R[\dim R] \simeq D^{\cdot}_R$
for the dualizing complex $D^{\cdot}_R$ (cf. \cite{rH}). Therefore
$\varprojlim \HH^i(R/\mathfrak a^n) \simeq \Hom(H^{\dim R
-i}_{\mathfrak a}(R), E)$ (cf. \ref{3.5}), which proves (d).
\end{proof}

It is a difficult problem to determine the cohomological dimension
$\cd (\mathfrak a, R).$ So the above result (d) in \ref{4.7}
illustrates the difficulty in order to calculate $\fgrade (
\mathfrak a, R).$ In the next result there is a generalization of
\ref{4.7} (d) for an arbitrary finitely generated module
$R$-module $M.$

\begin{theorem} \label{4.z} Let $(R, \mathfrak m)$ denote a local ring
with a dualizing complex $D^{\cdot}_R.$ Let $\mathfrak a$ denote
an ideal of $R.$ Then
\[
\fgrade (\mathfrak a, M) = \inf \{ i - \cd (\mathfrak a, K^i(M)) :
i = 0, \ldots , \dim M\}
\]
for a finitely generated $R$-module $M.$
\end{theorem}

\begin{proof} By the definition of the formal grade and Theorem
\ref{3.5} there is the equality
\[
\fgrade (\mathfrak a, M) = - \sup \{i \in \mathbb Z : H^i_{\mathfrak a}(\Hom (M, D^{\cdot}_R))
\not= 0 \}.
\]
Let $X$ denote an arbitrary complex of $R$-modules. Put
$s(X) = \sup \{i \in \mathbb Z : H^i(X) \not= 0\}.$ Let $\xx = x_1,
\ldots, x_r$ denote a system of elements of $R$ generating the ideal
$\mathfrak a.$ Let $\Cech$ denote the corresponding \v{C}ech complex.
Then
\[
H^i_{\mathfrak a}(\Hom (M,D^ {\cdot}_R)) \simeq H^i(\Cech \otimes
\Hom (M, D^{\cdot}_R))
\]
for all $i \in \mathbb Z$  (cf. \cite[Theorem 3.2]{pS2}). Therefore, it will
be enough to compute
\[
s(\Cech \otimes \Hom (M, D^{\cdot}_R)).
\]
Since $\Hom(M,D^{\cdot}_R)$
is a bounded complex with finitely generated cohomology modules and $\Cech$ is a bounded
complex of flat $R$-modules it follows that
\[
s(\Cech \otimes \Hom (M, D^{\cdot}_R)) =
\sup \{s(\Cech \otimes H^i(\Hom (M,D^ {\cdot}_R)) ) + i : i \in \mathbb Z\}
\]
(cf. \cite[Proposition 2.5]{hF1}). Because of $K^{-i}(M) = H^i(\Hom(M,D^{\cdot}_R)),
i \in \mathbb Z,$ it turns out that $s(\Cech \otimes H^i(\Hom (M,D^ {\cdot}_R)) = \cd (
\mathfrak a, K^{-i}(M)$ by the definition of the cohomological dimension. Whence the
claim is shown to be true.
\end{proof}

There is an expression of the cohomological
dimension in terms of the cohomological dimension of the minimal primes (cf. Corollary \ref{5.4}).
One might expect a similar result for the formal grade expressing
$\fgrade (\mathfrak a, M)$ in terms of the minimum of $\fgrade( \mathfrak a, R/\mathfrak p),$
where the minimum is taken over all $\mathfrak p \in \Min M$ or
$\mathfrak p \in \Ass M.$ This is not the case as the following example shows.

\begin{example} \label{5.5} Let $(R, \mathfrak m)$ denote a $d$-dimensional complete
local domain such that $H^i_{\mathfrak m}(R) = 0$ for all $i \not= 1, d, H^1_{\mathfrak m}(R)
\simeq k$ and $d \geq 4.$ Such rings exist. Let $D$ denote the global transform of $R.$
Then $D$ is a finitely generated
$R$-module with $H^i_{\mathfrak m}(D) = 0$ for all $i \not= d.$  Then $K(R) \simeq K(D)$
as easily seen. Now choose $\{x,y\}$ a $K(D)$-regular sequence and $\mathfrak a = (x,y)R.$
It follows that $\fgrade(\mathfrak a, D) = d- 2, \fgrade(\mathfrak a, R) = 1$ (cf. \ref{4.z}),
while $\Ass R = \Ass D = \{(0)\}.$
\end{example}

Moreover the example also shows that there are local rings such
that $\fgrade (\mathfrak a, R) \not= \dim R - \cd(\mathfrak a,
K(R)).$ But in any case there is the following bound for the
formal grade.

\begin{corollary} \label{4.a} Let $\mathfrak a$ be an ideal of the
local ring $(R, \mathfrak m).$ Then
\[
\fgrade (\mathfrak a, M) \leq \dim M - \cd(\mathfrak a, M)
\]
for a finitely generated $R$-module $M.$
\end{corollary}

\begin{proof} By Corollary \ref{5.4} there exists a prime ideal
$\mathfrak p \in \Ass_R M$ such that $\cd(\mathfrak a, M) = \cd (\mathfrak a,
R/\mathfrak p).$ Moreover, it follows that $\mathfrak p \in \Ass K^i(M)$
for a certain $0 \leq i \leq \dim M,$ (cf. Proposition \ref{2.y}). But
this implies $\cd(\mathfrak a, R/\mathfrak p) \leq \cd(\mathfrak a, K^i(M))$
as it is again a consequence of Corollary \ref{5.4}. By Theorem \ref{4.z} this
implies that
\[
\fgrade(\mathfrak a, M) \leq \dim M - \cd(\mathfrak a, K^i(M))
\leq \dim M - \cd(\mathfrak a, M),
\]
as required.
\end{proof}

Because $\height_M \mathfrak a \leq \cd(\mathfrak a, M)$ it
follows that the bound in Corollary \ref{4.a} is in fact an
improvement of the inequality $\fgrade(\mathfrak a, M) \leq \dim
M/\mathfrak aM$ (cf. Theorem \ref{4.5}).

Another difficulty about the formal grade is to characterize the
equality in \ref{4.7} (a). This has to do with a lack of
information about the $R$-module structure of $\varprojlim
\HH^i(M/\mathfrak a^nM), i \in \mathbb Z.$

\begin{theorem} \label{4.8} Let $M$ be a finitely generated
$R$-module. Then
\[
\dim \hat{R}/(\mathfrak a \hat{R},\mathfrak p) \geq \fgrade
(\mathfrak a, M)
\]
for all $\mathfrak p \in \Ass \hat{M}.$
\end{theorem}

\begin{proof} Without loss of generality one may assume that $R =
\hat{R}$ (cf. \ref{3.3}).  We proceed by induction on $t = \fgrade
(\mathfrak a, M).$ First consider the case of $t = 1.$ By our
assumption
\[
\Supp_R u_M(\mathfrak a) = \emptyset
\]
(cf. \ref{4.1}). But $\Supp_R u_M(\mathfrak a) = \cup_{\mathfrak p
\in T_{\mathfrak a}(M)}V(\mathfrak p)$ (cf. \ref{4.1}). This
implies that $\dim R/(\mathfrak a, \mathfrak p) \geq 1$ for all
$\mathfrak p \in \Ass_R M .$

Now let $t >1,$  i.e. in particular $\dim R/(\mathfrak a,
\mathfrak p) \geq 1$ for all $\mathfrak p \in \Ass_R M .$ By prime
avoidance arguments one may choose an element $x \in \mathfrak m$
which forms a parameter for all $R$-modules $R/(\mathfrak a,
\mathfrak p),$ where $\mathfrak p \in \Ass M.$

The long exact sequence
\[
\ldots \to \Hom(R_x,\varprojlim \HH^i(M/\mathfrak a^nM)) \to
\varprojlim \HH^i(M/\mathfrak a^nM) \to \varprojlim
\HH^i(M/(\mathfrak a,x)^nM) \to \ldots
\]
(cf.\ref{3.8}) provides that $\varprojlim \HH^i(M/(\mathfrak
a,x)^nM) = 0$ for all $i < t - 1.$ Therefore
\[
\dim R/(\mathfrak a, xR, \mathfrak p) \geq t -1 \text{ for all }
\mathfrak p \in \Ass M
\]
as a consequence of the the inductive hypothesis.

By the choice of $x \in \mathfrak m$ as a parameter for all
$R/(\mathfrak a, \mathfrak p), \mathfrak p \in \Ass M,$ this
proves that $\dim R/(\mathfrak a, \mathfrak p) \geq t$ for all
$\mathfrak p \in \Ass M.$ This completes the inductive step.
\end{proof}

In general the equality in Theorem \ref{4.8} does not hold. In
fact, this has to do with certain connectedness properties studied
in more detail in the next section.

\begin{example} \label{4.9} Let $R = k[|x_1,x_2,x_3,x_4|]$ denote
the formal power series ring in four variables over a field $k.$
Put $\mathfrak c = (x_1,x_2)R \cap (x_3,x_4)R.$ Then $\fgrade(
\mathfrak c, R) = 1$ (cf. Example \ref{5.2}), while $\dim
R/\mathfrak c = 2.$
\end{example}

We will continue here with another estimate of the formal grade
related to the cohomological dimension of ceratin associated prime
ideals.

\begin{theorem} \label{4.10} Let $(R,\mathfrak m)$ be a local ring.
Let $M$ denote a finitely generated $R$-module. Then
\[
\dim \hat{R}/\mathfrak p \geq \cd(\mathfrak a \hat{R},
\hat{R}/\mathfrak p) + \fgrade(\mathfrak a, M)
\]
for all $\mathfrak p \in \Ass_{\hat{R}}\hat{M}.$
\end{theorem}

\begin{proof} As mentioned above we may assume $R = \hat{R}$ as
follows by passing to the completion (cf. \ref{3.3}).  Now let
$\mathfrak p \in \Ass M$ be an associated prime ideal with $\dim
R/\mathfrak p = i$ for a certain $0 \leq i \leq \dim M.$ That is
\[
\mathfrak p \in (\Ass M)_i = (\Ass K^i(M))_i \; \mbox{and} \; \dim
K^i(M) = i
\]
(cf. \ref{2.y}). Moreover, it follows that $\Supp R/\mathfrak p
\subseteq \Supp K^i(M).$ Therefore (cf. \ref{5.3}) we see that
$\cd (\mathfrak a, R/\mathfrak p) \leq \cd (\mathfrak a, K^i(M)).$

By the assumption and the conclusion above it follows
\[
i - \cd (\mathfrak a, R/\mathfrak p) \geq i -\cd (\mathfrak a,
K^i(M)) \geq \fgrade(\mathfrak a, M)
\]
(cf. \ref{4.z}). Because of $i = \dim R/\mathfrak p$ this finishes
the proof.
\end{proof}

As $\hat{R}/\mathfrak p$ is a complete local domain it is a
catenary ring and therefore
\[
\dim \hat{R}/\mathfrak p = \dim \hat{R}/(\mathfrak a \hat{R},
\mathfrak p) + \height (\mathfrak a \hat{R}, \mathfrak
p)/\mathfrak p.
\]
Moreover $\height \height (\mathfrak a \hat{R}, \mathfrak
p)/\mathfrak p \leq \cd ( \mathfrak a \hat{R}, \hat{R}/\mathfrak
p).$ So, Theorem \ref{4.10} is in fact a sharpening of Theorem
\ref{4.8}.

\section{Connectedness properties}
\subsection{The Mayer-Vietoris sequence}
As it is well-known (cf. e.g. \cite[Section 19]{BS}, \cite{HH} and
\cite{pS2}) the Mayer-Vietoris sequence in local cohomology is an
important tool for connectedness phenomenons. Here we want to
continue with a variant of the Mayer-Vietoris sequence for formal
cohomology.

\begin{theorem} \label{5.1} Let $\mathfrak a, \mathfrak b$ two
ideals of a local ring $(R, \mathfrak m).$ For a finitely
generated $R$-module $M$ there is the long exact sequence
\begin{multline*}
\ldots \to \varprojlim \HH^i(M/(\mathfrak a \cap \mathfrak b)^nM)
\to \varprojlim \HH^i(M/\mathfrak a^nM) \oplus \varprojlim
\HH^i(M/\mathfrak b^nM) \to
\\ \to \varprojlim \HH^i(M/(\mathfrak a,
\mathfrak b)^nM) \to \ldots,
\end{multline*}
where $i \in \mathbb Z.$
\end{theorem}

\begin{proof} Let $n \in \mathbb Z$ denote an integer. Then there
is the the following natural exact sequence
\[
0 \to M/(\mathfrak a^nM \cap \mathfrak b^n M) \to M/\mathfrak a^nM
\oplus M/\mathfrak b^nM \to M/(\mathfrak a^n, \mathfrak b^n)M \to
0.
\]
Now the long exact local cohomology sequence provides by passing
to the projective limit the following long exact cohomology
sequence
\begin{multline*}
\ldots \to \varprojlim \HH^i(M/(\mathfrak a^n M\cap \mathfrak
b^nM)) \to \varprojlim \HH^i(M/\mathfrak a^nM) \oplus \varprojlim
\HH^i(M/\mathfrak b^nM) \to \\ \to \varprojlim \HH^i(M/(\mathfrak
a^n, \mathfrak b^n)M) \to \ldots.
\end{multline*}
Notice that the projective limit on projective systems of Artinian
modules is exact.

Now we observe that the $(\mathfrak a, \mathfrak b)$-adic
filtration is equivalent to the filtration $\{(\mathfrak a^n,
\mathfrak b^n)M\}_{n \in \mathbb N}.$ In order to finish the proof
we have to show that the $(\mathfrak a \cap \mathfrak b)$-adic
filtration on $M$ is equivalent to the filtration $\{(\mathfrak
a^n \cap \mathfrak b^n)M)\}_{n \in \mathbb N}$ (cf. \ref{3.6a}).

To this end first note that $(\mathfrak a \mathfrak b)^nM
\subseteq (\mathfrak a^n \cap \mathfrak b^n)M \subseteq \mathfrak
a^nM \cap \mathfrak b^n M$ for all $n \in \mathbb N.$ Let $m \in
\mathbb N$ denote a given integer. By the Artin-Rees Lemma (cf.
\cite[Ch. III, \S 3, Cor. 1]{nB}) there exists an $k \in \mathbb
N$ such that $\mathfrak a^n N\cap \mathfrak b^m N\subseteq
\mathfrak a^{n-k}\mathfrak b^mN$ for all $n \geq k.$ Since the
$\mathfrak a \mathfrak b$-adic and the $\mathfrak a \cap \mathfrak
b$-adic topology on $M$ are equivalent this finishes the proof.
\end{proof}

The above result (cf. \ref{5.1}) provides an example related to
the supports of formal cohomology.

\begin{example} \label{5.2} Let $k$ be a field. Let $R =
k[|x_1,x_2,x_3,x_4|]$ denote the formal power series ring in four
variables over $k.$ Put $\mathfrak a = (x_1,x_2)R$ and $\mathfrak
b = (x_3,x_4)R.$ Then the Mayer-Vietoris sequence provides the
following two isomorphisms
\[
R \simeq \varprojlim \HH^1(R/(\mathfrak a \cap \mathfrak b)^n)\,
\mbox{and} \, \varprojlim \HH^2(R/(\mathfrak a \cap \mathfrak
b)^n) \simeq \varprojlim \HH^2(R/\mathfrak a^n) \oplus \varprojlim
\HH^2(R/\mathfrak b^n).
\]
To this end remark that $(\mathfrak a, \mathfrak b)$ is the
maximal ideal of the complete local ring $R.$ Therefore $\Supp
\HH^1(R/(\mathfrak a \cap \mathfrak b)^n) = \Spec R,$ while  $\dim
R/\mathfrak a \cap \mathfrak b = 2.$
\end{example}

Note that the example was introduced by Hartshorne (cf.
\cite{rH1}). In the following we want to extend these
considerations to a more subtle investigation.

\subsection{On the connectedness}
Next let us summarize a few technical preparations for the
connectedness results. Let $(R, \mathfrak m)$ denote a local ring.

\begin{lemma} \label{5.6} Let $M$ be a finitely generated $R$-module. Let $\mathfrak a, \mathfrak b$
denote two ideals of $R.$ Suppose that $\varprojlim
\HH^1(M/(\mathfrak a \cap \mathfrak b)^ nM) = 0.$ Then
$T_{\mathfrak a \hat R}(\hat{M}) \cup T_{\mathfrak b \hat
R}(\hat{M}) = T_{(\mathfrak a, \mathfrak b)\hat R}(\hat{M}).$
\end{lemma}

\begin{proof} First remember that we may assume that $(R, \mathfrak m)$ is a complete local ring
(cf. \ref{3.3}). With the notation introduced in Section 4.1
it is clear that the left hand side of the
statement is contained in the right hand side.

In order to prove the reverse containment relation the
Mayer-Vietoris sequence (cf. \ref{5.1}) provides an epimorphism
\[
u_M(\mathfrak a)\oplus u_M(\mathfrak b) \to u_M(\mathfrak a, \mathfrak b) \to 0
\]
(use Lemma \ref{4.1}). Now let $\mathfrak p \in \Ass u_M(\mathfrak a, \mathfrak b),$
i.e. $\mathfrak p \in \Ass M$ and $\dim R/(\mathfrak p, \mathfrak a, \mathfrak b) = 0.$ In particular
it follows that $\mathfrak p \in \Supp u_M(\mathfrak a, \mathfrak b)$ and therefore
$\mathfrak p \in \Supp u_M(\mathfrak a) \oplus u_M(\mathfrak b).$ Without loss of generality we may
conclude that $\mathfrak p \in \Supp u_M(\mathfrak a).$ So there exists a prime ideal
$\mathfrak q \in \Ass u_M(\mathfrak a)$ with $\mathfrak q \subseteq \mathfrak p.$ Whence
$\mathfrak q \in \Ass M$ and $\dim R/(\mathfrak q, \mathfrak a) = 0$ (cf. Lemma \ref{4.1}).
Because of $\mathfrak p \in \Ass M$ and $\mathfrak q \subseteq \mathfrak p$ this implies
$\mathfrak p \in \Ass u_M(\mathfrak a),$ which finishes the proof.
\end{proof}

As another consequence of the Mayer-Vietoris sequence there is the
following connectedness result. To this end an $R$-module $M$ is
called indecomposable whenever $M = N_1 \oplus N_2$ implies either
$M = N_1$ and $N_2 = 0$ or $N_1 = 0$ and $M = N_2.$

\begin{lemma} \label{5.7} Let $\hat{M}$ denote an indecomposable $\hat{R}$-module.
Suppose that $\fgrade (\mathfrak a, M) \geq 2$ for an ideal
$\mathfrak a$ of $R.$ Then $\Supp_{\hat R} {\hat M}/\mathfrak a
{\hat M} \setminus \{\hat{\mathfrak m}\}$ is connected.
\end{lemma}

\begin{proof} Because of $\fgrade (\mathfrak a, M) = \fgrade (\mathfrak a \hat R, \hat M)$
(cf. \ref{3.3}) we may assume that $R$ is a complete local ring.
Now suppose that $\Supp M/\mathfrak a M \setminus \{\mathfrak m\}$ is disconnected. Then
there are two ideals $\mathfrak b, \mathfrak c$ of $R$ satisfying the following properties
\begin{itemize}
\item[1.] $\Rad (\mathfrak a, \Ann M) = \Rad( \mathfrak b \cap \mathfrak c),$
\item[2.] $(\mathfrak b, \mathfrak c)$ is an $\mathfrak m$-primary ideal, and
\item[3.] neither $\mathfrak b$ nor $\mathfrak c$ is an $\mathfrak m$-primary ideal.
\end{itemize}
Then the Mayer-Vietoris sequence (cf. \ref{5.1}) provides an
isomorphism
\[
\varprojlim \HH^0(M/\mathfrak b^nM) \oplus \varprojlim
\HH^0(M/\mathfrak c^nM) \simeq \varprojlim \HH^0(M/(\mathfrak b,
\mathfrak c)^nM).
\]
But $(\mathfrak b, \mathfrak c)$ is an $\mathfrak m$-primary ideal and therefore
$\varprojlim \HH^0(M/(\mathfrak b,\mathfrak c)^nM) \simeq M$ (cf. \ref{4.4}).
By the indecomposability of $M$ it follows -- say --
\[
\varprojlim \HH^0(M/\mathfrak b^nM) \simeq M \; \mbox{ and } \;
\varprojlim \HH^0(M/\mathfrak c^nM) = 0.
\]
Therefore, by \ref{4.1} it turns out that $\dim R/(\mathfrak p, \mathfrak b) = 0$ for
all $\mathfrak p \in \Ass M.$ Because of
\[
\mathfrak m = \cap_{\mathfrak p \in \Ass M} \Rad (\mathfrak b,
\mathfrak p) = \Rad \mathfrak b
\]
it yields that $\mathfrak b$ is an $\mathfrak m$-primary
ideal. This is a contradiction.
\end{proof}

One might observe that for the proof of $\Rad \mathfrak b =
\mathfrak m$ it will be enough to consider only the minimal prime
ideals $\mathfrak p \in  \Ass M.$ This is a corner stone for a
generalization in the next subsection.

The indecomposibility of $M$ in \ref{5.7} is essential as the
following example shows.

\begin{example} \label{5.8} With the notation of Example \ref{5.2} put
$M = R/\mathfrak a \oplus R/\mathfrak b.$ Let $\mathfrak c =
\mathfrak a \cap \mathfrak b.$ Then $\fgrade (\mathfrak c, M) =
\depth M = 2,$ while $\Supp M/\mathfrak cM \setminus \{\mathfrak
m\}$ is not connected. Recall that $\mathfrak c = \Ann M.$
\end{example}
We apply the previous Lemma in order to derive a corresponding
connectedness result related to the cohomological dimension. To
this end we introduce the notion $\Assh M = \{ \mathfrak p \in
\Ass M : \dim R/\mathfrak p = \dim M\}$ for a finitely generated
$R$-module.

\begin{theorem} \label{5.7a} Let $(R, \mathfrak m)$ denote a local
ring. Let $\mathfrak a$ be an ideal of $R.$ Suppose that
\begin{itemize}
\item[(a)] $\Ass \hat R = \Assh \hat R,$

\item[(b)] $H^{\dim R}_{\mathfrak m}(R)$  is indecomposable,

\item[(c)] $\cd(\mathfrak a, R) \leq \dim R -2.$
\end{itemize}
Then $V(\mathfrak a \hat R) \setminus V(\mathfrak m \hat R)$ is
connected.
\end{theorem}

\begin{proof} Because of $\cd (\mathfrak a, R) = \cd (\mathfrak a \hat
{R} , \hat R)$ one may assume that $R$  possesses a dualizing
complex (cf. \ref{3.3}). Observe that $H^d_{\mathfrak m}(R) \simeq
H^d_{\hat {\mathfrak m}}(\hat R), d = \dim R = \dim \hat R.$

Let $\mathbb Q(R)$ denote the total ring of quotients of $R.$ Then
there exists a birational extension ring $R \subset S \subset
\mathbb Q(R)$ such that $S$ is a finitely generated $R$-module and
satisfies the condition $S_2$ (cf. \cite[5.3]{pS4}). To this end
we have to use (a). Whence it follows that
\[
\cd (\mathfrak a, K^i(S)) \leq \dim K^i(S) \leq i-2
\]
for all $0 \leq i < \dim S = d$ (cf. \ref{2.y}).
Moreover, the short exact sequence
\[
0 \to R \to S \to S/R \to 0
\]
provides the vanishing $H^i_{\mathfrak a}(S) = 0$ for all $i > d
-2.$ To this end observe that $\dim S/R \leq d -2$ (cf.
\cite[5.3]{pS4}) and that $H^i_{\mathfrak a}(R) = 0$ for $i >
d-2.$ Therefore $\cd (\mathfrak a, S) \leq \dim S -2.$ Since
$\Supp S = \Supp K(S)$ we obtain $\cd (\mathfrak a, S) = \cd
(\mathfrak a, K(S))$ (cf. \ref{5.4}). But then it follows that
\[
\fgrade (\mathfrak a, S) = \min \{i - \cd (\mathfrak a, K^i(S)) :
i = 0, \ldots, \dim S\} \geq 2
\]
(cf. \ref{4.z}). In order to apply \ref{5.7} we show that $S$ as an
$R$-module is indecomposable.

Assume the contrary, i.e.  $S \simeq S_1 \oplus S_2$ for two
non-zero $R$-modules $S_i, i= 1,2.$ Clearly $\dim S_i = d, i =
1,2.$ This follows since $S$ has the property that $\dim
S/\mathfrak p = \dim S$ for all $\mathfrak p \in \Supp_R S$ (cf.
\cite{pS4}).

By considering the local cohomology modules we see that
\[
\HH^d(R) \simeq \HH^d(S) \simeq \\H^d(S_1) \oplus \HH^d(S_2), \;
\mbox{ and } \; \HH^d(S_i) \not= 0, \; i = 1,2.
\]
Notice that $\dim S/R \leq d-2.$  Because $H^d_{\mathfrak m}(R)$
is supposed to be indecomposable by condition (b) this is a
contradiction.

So, the previous result (cf. \ref{5.7}) finally
implies that
\[
\Supp_R S/\mathfrak a S \setminus V(\mathfrak m) = V(\mathfrak a)
\setminus V(\mathfrak m)
\]
is connected. To this end remember that $\Supp_R S = \Spec R.$
\end{proof}

We note that Theorem \ref{5.7a} extends \cite[2.27]{pS1}, where
the condition $S_2$ is assumed for $R$ in order to derive the
connectedness property.
Note that the indecomposibility of 
$H^{\dim R}_{\mathfrak m}(R)$ was studied by Hochster and Huneke 
(cf. \cite[Theorem 4.1]{HH}).

\subsection{The connectedness dimension}
Next let us summarize a few technical preparations for further
connectedness results. Let $(R, \mathfrak m)$ denote a local ring.

\begin{definition} \label{5.9}
For an $R$-module$M$ define
\[
c(M) = \min \{\dim R/\mathfrak c : V(\mathfrak c)
\subseteq \Supp M \mbox{ and } \Supp M
\setminus V(\mathfrak c) \mbox{ is disconnected}\}.
\]
\end{definition}

We refer to \cite[Section 19]{BS} for more details about the
definition. Here we notice that $c(M) \leq \dim M$ with equality
provided $\Supp M$ is irreducible. Moreover $c(M) \geq 0.$

Now let $M$ be a finitely generated $R$-module. Let $\mathfrak
p_1, \ldots, \mathfrak p_r$ denote the distinct minimal prime
ideals of $\Supp M = V(\Ann_R M).$

Let $\mathcal S (r)$ denote the set of all ordered pairs $(A, B)$
of non-empty subsets of $\{1,\ldots,r\}$ such that $A\cup B =
\{1,\ldots,r\}.$

\begin{lemma} \label{5.10} Let $M$ be a finitely generated
$R$-module. Then
\[
c(M) = \min \{ \dim R/((\cap_{i \in A}\mathfrak p_i), (\cap_{j
\in B}\mathfrak p_j)) : (A,B) \in \mathcal S(r) \}.
\]
\end{lemma}

\begin{proof} The result is a module theoretic version of
\cite[19.1.15]{BS}. For the details of the proof we refer to
\cite[19.1.15 and 19.2.5]{BS}. To this end observe that $\Supp M =
V(\Ann M).$
\end{proof}

Next we want to continue with a the behavior of the connectedness dimension
by a generic hyperplane section. To be more precise:

\begin{lemma} \label{5.11} Let $M$ denote a finitely generated $R$-module
with $c(M) > 0.$ Then there exists an element $x \in \mathfrak m$
such that $c(M) \geq c(M/xM) +1.$
\end{lemma}

\begin{proof} Let $\mathfrak p_1, \ldots, \mathfrak p_r$ denote the distinct
minimal prime ideals of $V(\Ann_R M).$ Then $c(M) = \dim
R/\mathfrak c > 0$ for an ideal $\mathfrak c = ((\cap_{i \in
A}\mathfrak p_i), (\cap_{j \in B}\mathfrak p_j))$ with a certain
pair $(A,B) \in \mathcal S(r)$ (cf. Lemma \ref{5.10}). Now choose
$x \in \mathfrak m$ as a parameter of $R/\mathfrak c,$ i.e. $c(M)
-1 = \dim R/(xR, \mathfrak c).$

Next observe that $V(x,\mathfrak a \cap \mathfrak b) = V(x,
\mathfrak a) \cup V(x, \mathfrak b) = V((x,\mathfrak a)\cap
(x,\mathfrak b))$ for two ideals $\mathfrak a, \mathfrak b$ of
$R.$ Then there are the following equalities for the radical
ideals
\begin{equation*}
\begin{aligned}
\Rad (xR, \mathfrak c) & = \Rad (\cap_{i \in A}(\mathfrak p_i,
xR), \cap_{j \in B}(\mathfrak p_j, xR)) \\
  & =  \Rad (\cap_{i \in
A}\Rad(\mathfrak p_i, xR), \cap_{j \in B}(\Rad (\mathfrak p_j,
xR)))
\end{aligned}
\end{equation*}
as easily seen. Let $\mathfrak P_1, \ldots, \mathfrak P_s$ denote
the distinct minimal prime ideals of $V(xR,\Ann_R M).$ By easy
computations it follows that
\[
V(\Ann_R M, xR) = V(\cap_{i=1}^r \mathfrak p_i) \cap V(xR) =
V(\cap_{i=1}^r (\mathfrak p_i,xR)).
\]
Whence, the set of prime ideals $\mathfrak P_1, \ldots, \mathfrak
P_s$ coincides with the set of minimal prime ideals of the ideal
$\cap_{i=1}^r (\mathfrak p_i, xR)$ and
\[
\cap_{i=1}^s \mathfrak P_i = \Rad (\cap_{i=1}^r (\mathfrak
p_i,xR)) = \Rad (\Ann_R M, xR).
\]
By avoiding redundant components in  $\cap_{i \in A}\Rad(\mathfrak
p_i, xR)$ and $\cap_{j \in B}\Rad (\mathfrak p_j, xR)$ resp. we
derive a representation
\[
\cap_{i \in A}\Rad((\mathfrak p_i, xR)) = \cap_{i \in
\tilde{A}} \mathfrak P_i \quad \mbox{ and } \quad \cap_{j \in
B}\Rad((\mathfrak p_j, xR)) = \cap_{j \in \tilde{B}} \mathfrak
P_j
\]
for an ordered pair $(\tilde{A}, \tilde{B}) \in \mathcal S(s).$
This means that
\[
c(M) -1 = \dim R/(xR, \mathfrak c) = \dim R/(\cap_{i \in
\tilde{A}} \mathfrak P_i, \cap_{j \in \tilde{B}} \mathfrak P_j)
\geq c(M/xM),
\]
as required. Note that the dimension does not change by passing to the radical.
\end{proof}

As a consequence of the Lemmas \ref{5.11} and \ref{2.x} one has
the following result, relating the connectedness dimension of
$R/\mathfrak a$ and the cohomological dimension.

\begin{corollary} \label{3.a} Let $\mathfrak a$ be an ideal of a
local ring $(R, \mathfrak m).$ Suppose that $\HH^d(R)$ is
indecomposable and $\Ass \hat{R} = \Assh \hat{R}.$ Then
$c(\hat{R}/\mathfrak a \hat{R}) \geq \dim R - \cd(\mathfrak a, R)
-1.$
\end{corollary}

\begin{proof} First note that we may assume that $R = \hat{R},$
that is $R$ is complete (cf. \ref{3.3}).  For the proof we proceed
by an induction on $c(R/\mathfrak a).$ In the case of
$c(R/\mathfrak a) = 0$ the result is a consequence of \ref{5.7a}.
So assume that $c(R/\mathfrak a) > 0.$ Then there exists an
element $x \in \mathfrak m$ such that $c(R/\mathfrak a) \geq
c(R/(\mathfrak a, xR)) +1$ (cf. \ref{5.11}). By the inductive
hypothesis
\[
c(R/\mathfrak a)-1 \geq c(R/(\mathfrak a, xR)) \geq \dim R - \cd
((\mathfrak a, xR),R) -1.
\]
On the other hand $\cd ((\mathfrak a, xR),R) \leq \cd (\mathfrak
a, R) +1$ (cf. \ref{2.x}). Now this completes the inductive step
by putting together  these inequalities.
\end{proof}

In their paper \cite[Theorem 3.4]{DNT} the authors claimed the
validity of \ref{3.a} without the condition that $H^d_{\mathfrak
m}(R)$ is indecomposable. This is not correct as follows by
Example \ref{5.2}. To this end let  $\mathfrak c = \mathfrak a
\cap \mathfrak b.$ Then $\cd (\mathfrak c, R/\mathfrak c) = 0,
\dim R/\mathfrak c = 2, c(R/\mathfrak c) = 0.$ Moreover
$\HH^2(R/\mathfrak c) \simeq \HH^2(R/\mathfrak a) \oplus
\HH^2(R/\mathfrak b)$ and both of the direct summands do not
vanish.

\subsection{Formal cohomology and connectedness}
In this subsection we relate the vanishing the formal cohomology
to the connectedness properties.

\begin{theorem} \label{5.12} Let $\mathfrak a$ denote an ideal of a local
ring $(R, \mathfrak m).$ Let $M$ be a finitely generated
$R$-module. Then $c(\hat{R}/(\mathfrak a \hat{R}, \mathfrak p))
\geq \fgrade (\mathfrak a, M)  - 1$ for all $\mathfrak p \in
\Ass_{\hat R} \hat{M}.$
\end{theorem}

\begin{proof}
First of all we note that Corollary \ref{3.a} applied to
$\mathfrak a$ in $\hat{R}/\mathfrak p, \mathfrak p \in
\Ass_{\hat{R}}\hat{M},$ provides the following inequality
\[
c(\hat{R}/(\mathfrak a \hat{R}, \mathfrak p)) \geq \dim
\hat{R}/\mathfrak p - \cd(\mathfrak a \hat{R}, \hat{R}/\mathfrak
p) -1.
\]
To this end we have to prove that $\HH^i(\hat{R}/\mathfrak p), i =
\dim \hat{R}/\mathfrak p,$ is indecomposable. By local duality it
will be enough to prove that the canonical module
$K(\hat{R}/\mathfrak p)$ is an indecomposable $\hat{R}/\mathfrak
p$-module. Since $\hat{R}/\mathfrak p$ is a domain and since
$K(\hat{R}/\mathfrak p)$ is a torsion-free $\hat{R}/\mathfrak
p$-module of rank 1, it is indecomposable.

 On the other hand (cf. \ref{4.8}) it follows that
\[
\dim \hat{R}/\mathfrak p - \cd(\mathfrak a \hat{R},
\hat{R}/\mathfrak p) \geq \fgrade(\mathfrak a, M).
\]
Putting together both of the estimates the desired inequality is
shown to be true.
\end{proof}

As a particular case of Theorem \ref{5.12} there is the following corollary.

\begin{corollary} \label{5.13} Let $\mathfrak a$ denote an ideal of a local
ring $(R, \mathfrak m).$ Let $M$ be a finitely generated
$R$-module. Suppose that $\varprojlim \HH^i(M/\mathfrak a^nM) = 0$
for $i \leq 1.$

Then $V(\mathfrak a \hat{R}, \mathfrak p) \setminus V(\hat{\mathfrak m})$
is connected for all $\mathfrak p \in \Ass_{\hat R} \hat{M}.$
\end{corollary}

\begin{proof} As follows by the definitions the claim is a particular case of \ref{5.12}.
To this end recall that $\fgrade( \mathfrak a, M) \geq 2.$
\end{proof}

It is noteworthy to remark that the converse of the previous results are not true.

\begin{example} \label{5.14} With the notion of \ref{5.2} put $M = R/\mathfrak c,
\mathfrak c = \mathfrak a \cap \mathfrak b.$ Then $V(\mathfrak c, \mathfrak p) \setminus
V(\mathfrak m)$ is connected for all $\mathfrak p \in \Ass M,$ while
\[
\varprojlim \HH^1(M/\mathfrak c^n M) \simeq \HH^1(M) \simeq R/\mathfrak m,
\]
as it is easily seen.
\end{example}

As further application of the results of this and the previous
subsection there is another estimate of the formal grade, more in
the sense of Theorem \ref{5.7}.

\begin{corollary} \label{5.15} Let $M$ denote a finitely generated
$R$-module, where $(R, \mathfrak m)$ is a local ring. Suppose that
\begin{itemize}
\item[(a)] $\Ass_{\hat R}\hat{M} = \Assh_{\hat R}\hat{M}$ and
\item[(b)] $\HH^d(R/\Ann_R M), d = \dim M,$ is indecomposable.
\end{itemize}
Then $c(\hat{M}/\mathfrak a \hat{M}) \geq \fgrade(\mathfrak a, M)
-1.$
\end{corollary}

\begin{proof} Without loss of generality we may assume that $R$ is
a complete local ring (cf. \ref{3.3}). Moreover, by the definition
it follows that $c(M/\mathfrak aM) = c(R/(\mathfrak a, \Ann_R
M)).$ The assumption (a) implies that $\Ass R/\Ann_R M = \Assh
R/\Ann_R M.$ Because $\HH^d(R/\Ann_R M)$ is indecomposable we may
apply \ref{5.7}, so that
\[
c(R/(\mathfrak a,\Ann_R M)) \geq \dim R/\Ann_R M - \cd(\mathfrak
a, R/\Ann_R M) -1.
\]
But now $\dim M = \dim R/\Ann_R M.$ Furthermore $\cd(\mathfrak
a,M) = \cd(\mathfrak a, R/\Ann_R M)$ (cf. \ref{5.4}). Because of
$\dim M - \cd(\mathfrak a, M) \geq \fgrade(\mathfrak a, M)$ (cf.
\ref{4.a}) this finishes the proof.
\end{proof}

\end{document}